\documentclass[12pt]{article}
\usepackage{amssymb, amsmath, amsthm, amscd}
\usepackage[pdftex]{graphics}
\usepackage[utf8]{inputenc}
\usepackage[all,cmtip]{xy}
\usepackage{bbm}
\usepackage{enumitem}
\usepackage{setspace}
\usepackage{mathdots}

\usepackage[mathscr]{euscript}
\usepackage[colorlinks=true]{hyperref}
\hypersetup{colorlinks   = true}
\hypersetup{linkcolor=blue}
\usepackage{tikz}
\usetikzlibrary{cd}

\begin{document}

\mathchardef\mhyphen="2D
\newtheorem{The}{Theorem}[section]
\newtheorem{Lem}[The]{Lemma}
\newtheorem{Prop}[The]{Proposition}
\newtheorem{Cor}[The]{Corollary}
\newtheorem{Rem}[The]{Remark}
\newtheorem{Obs}[The]{Observation}
\newtheorem{SConj}[The]{Standard Conjecture}
\newtheorem{Titre}[The]{\!\!\!\! }
\newtheorem{Conj}[The]{Conjecture}
\newtheorem{Question}[The]{Question}
\newtheorem{Prob}[The]{Problem}
\newtheorem{Def}[The]{Definition}
\newtheorem{Not}[The]{Notation}
\newtheorem{Claim}[The]{Claim}
\newtheorem{Conc}[The]{Conclusion}
\newtheorem{Ex}[The]{Example}
\newtheorem{Fact}[The]{Fact}
\newtheorem{Formula}[The]{Formula}
\newtheorem{Formulae}[The]{Formulae}
\newtheorem{The-Def}[The]{Theorem and Definition}
\newtheorem{Prop-Def}[The]{Proposition and Definition}
\newtheorem{Lem-Def}[The]{Lemma and Definition}
\newtheorem{Cor-Def}[The]{Corollary and Definition}
\newtheorem{Conc-Def}[The]{Conclusion and Definition}
\newtheorem{Terminology}[The]{Note on terminology}
\newcommand{\C}{\mathbb{C}}
\newcommand{\R}{\mathbb{R}}
\newcommand{\N}{\mathbb{N}}
\newcommand{\Z}{\mathbb{Z}}
\newcommand{\Q}{\mathbb{Q}}
\newcommand{\Proj}{\mathbb{P}}
\newcommand{\Rc}{\mathcal{R}}
\newcommand{\Oc}{\mathcal{O}}
\newcommand{\Vc}{\mathcal{V}}
\newcommand{\Id}{\operatorname{Id}}
\newcommand{\pr}{\operatorname{pr}}
\newcommand{\rk}{\operatorname{rk}}
\newcommand{\del}{\partial}
\newcommand{\delbar}{\bar{\partial}}
\newcommand{\Cdot}{{\raisebox{-0.7ex}[0pt][0pt]{\scalebox{2.0}{$\cdot$}}}}
\newcommand\nilm{\Gamma\backslash G}
\newcommand\frg{{\mathfrak g}}
\newcommand{\fg}{\mathfrak g}
\newcommand{\Oh}{\mathcal{O}}
\newcommand{\Kur}{\operatorname{Kur}}
\newcommand\gc{\frg_\mathbb{C}}
\newcommand\jonas[1]{{\textcolor{green}{#1}}}
\newcommand\luis[1]{{\textcolor{red}{#1}}}
\newcommand\dan[1]{{\textcolor{blue}{#1}}}

\begin{center}

{\Large\bf Balanced Hyperbolic and Divisorially Hyperbolic Compact Complex Manifolds}

\end{center}

\begin{center}

{\large Samir Marouani and Dan Popovici}

\end{center}

\vspace{1ex}

\noindent{\small{\bf Abstract.} We introduce two notions of hyperbolicity for not necessarily K\"ahler $n$-dimensional compact complex manifolds $X$. The first, called {\it balanced hyperbolicity}, generalises Gromov's K\"ahler hyperbolicity by means of Gauduchon's balanced metrics. The second, called {\it divisorial hyperbolicity}, generalises the Brody hyperbolicity by ruling out the existence of non-degenerate holomorphic maps from $\C^{n-1}$ to $X$ that have what we term a subexponential growth. Our main result in the first part of the paper asserts that every balanced hyperbolic $X$ is also divisorially hyperbolic. We provide a certain number of examples and counter-examples and discuss various properties of these manifolds. In the second part of the paper, we introduce the notions of {\it divisorially K\"ahler} and {\it divisorially nef} real De Rham cohomology classes of degree $2$ and study their properties. They also apply to $C^\infty$, not necessarily holomorphic, complex line bundles and are expected to be implied in certain cases by the hyperbolicity properties introduced in the first part of the work. While motivated by the observation of hyperbolicity properties of certain non-K\"ahler manifolds, all these four new notions seem to have a role to play even in the K\"ahler and the projective settings.}

\vspace{2ex}

\section{Introduction}\label{section:Introduction} We propose a hyperbolicity theory in which curves are replaced by divisors and the bidegree $(1,\,1)$ is replaced by the bidegree $(n-1,\,n-1)$ on $n$-dimensional compact complex manifolds. The notions we introduce are weaker, hence more inclusive, than their classical counterparts. In particular, the setting need not be projective or even K\"ahler. Recall that, by contrast, the classical notions of K\"ahler and Brody/Kobayashi hyperbolicity cannot occur in the context of compact manifolds, at least conjecturally for the latter, outside the projective context. Our motivation stems from the existence of many interesting examples of non-K\"ahler compact complex manifolds that display hyperbolicity features in the generalised sense that we now set out to explain.

\vspace{2ex}

{\bf (I)}\, Let $X$ be a compact complex manifold with $\mbox{dim}_\C X=n\geq 2$. Recall that

\vspace{2ex}

(a)\, $X$ is said to be {\it K\"ahler hyperbolic} in the sense of Gromov (see [Gro91]) if there exists a K\"ahler metric $\omega$ on $X$ whose lift $\widetilde\omega$ to the universal covering space $\widetilde{X}$ of $X$ is $d$-exact with an $\widetilde\omega$-bounded $d$-potential on $\widetilde{X}$.

Meanwhile, it is well known that balanced metrics (i.e. Hermitian metrics $\omega$ on $X$ such that $d\omega^{n-1}=0$) may exist on certain compact complex $n$-dimensional manifolds $X$ such that $\omega^{n-1}$ is even $d$-exact on $X$. These manifolds are called {\it degenerate balanced}. There is no analogue of this phenomenon in the K\"ahler setting. Building on this fact, we propose in Definition \ref{Def:bal-hyperbolic} the balanced analogue of Gromov's K\"ahler hyperbolicity by requiring the existence of a balanced metric $\omega$ on $X$ such that $\widetilde\omega^{n-1}$ is $d$-exact with an $\widetilde\omega$-bounded $d$-potential on $\widetilde{X}$. We call any compact manifold $X$ admitting such a metric $\omega$ a {\bf balanced hyperbolic} manifold. We immediately get our first examples: any degenerate balanced manifold $X$ is automatically balanced hyperbolic.

\vspace{2ex}

(b)\, $X$ is said to be {\it Kobayashi hyperbolic} (see e.g. [Kob70]) if the Kobayashi pseudo-distance on $X$ is actually a distance. By Brody's Theorem 4.1. in [Bro78], this is equivalent to the non-existence of entire curves in $X$, namely the non-existence of non-constant holomorphic maps $f:\C\longrightarrow X$. This latter property has come to be known as the {\it Brody hyperbolicity} of $X$. Thus, a compact manifold $X$ is Brody hyperbolic if and only if it is Kobayashi hyperbolic. (The equivalence is known to fail when $X$ is non-compact.) In Definition \ref{Def:div-hyperbolic}, we propose the $(n-1)$-dimensional analogue of the Brody hyperbolicity, that we call {\bf divisorial hyperbolicity}. However, due to the absence of a higher-dimensional analogue of Brody's Reparametrisation Lemma [Bro78, Lemma 2.1.], there is a surprising twist: to ensure that the {\it balanced hyperbolicity} of $X$ implies its {\it divisorial hyperbolicity}, it does not suffice to rule out the existence of non-degenerate (at some point) holomorphic maps $f:\C^{n-1}\longrightarrow X$ in Definition \ref{Def:div-hyperbolic} when $\mbox{dim}_\C X=n\geq 3$, but the non-existence has to be confined to such maps of {\it subexponential growth} in the sense of Definition \ref{Def:subexp-growth}.

This restriction is also warranted by the trivial existence of a non-degenerate (at some point) holomorphic map $f:\C^{n-1}\longrightarrow X$ whenever $X=G/\Gamma$ is the quotient of an $n$-dimensional {\it complex} Lie group $G$ by a discrete co-compact subgroup $\Gamma$. Indeed, such a map is obtained by composing two obvious maps with the {\it exponential} $\exp:\fg\longrightarrow G$ from the Lie algebra $\fg$ of $G$. However, we will see that some of these quotients $X=G/\Gamma$ deserve to be called {\it divisorially hyperbolic}.

\vspace{2ex}

To put our {\it balanced hyperbolicity} and {\it divisorial hyperbolicity} into perspective, we sum up below the relations among the various hyperbolicity notions mentioned above.

\begin{The}\label{The:introd_hyperbolicity-implications} Let $X$ be a compact complex manifold. The following implications hold:

\vspace{3ex}  

\hspace{12ex} $\begin{array}{lll} X \hspace{1ex} \mbox{is {\bf K\"ahler hyperbolic}} & \implies & X \hspace{1ex} \mbox{is {\bf Kobayashi/Brody hyperbolic}} \\
 \rotatebox{-90}{$\implies$} &  & \rotatebox{-90}{$\implies$} \\
 X \hspace{1ex} \mbox{is {\bf balanced hyperbolic}} & \implies & X \hspace{1ex} \mbox{is {\bf divisorially hyperbolic}} \\
 \rotatebox{90}{$\implies$} &  &   \\
 X \hspace{1ex} \mbox{is {\bf degenerate balanced}} & &\end{array}$ 

\end{The}

\vspace{2ex}

The vertical implications in Theorem \ref{The:introd_hyperbolicity-implications} are obvious, while the top horizontal implication has been known since [Gro91, 0.3.B.]. (See also [CY18, Theorem 4.1.] for a proof.) Our main result of $\S.$\ref{section:bal-div_hyperbolicity} is Theorem \ref{The:bal-div-hyperbolic_implication} proving the bottom horizontal implication of Theorem \ref{The:introd_hyperbolicity-implications}.

\vspace{2ex}

We now collect a few results from $\S.$\ref{subsection:bal-hyp_examp}. They contain various examples and counter-examples illustrating our {\it balanced hyperbolicity} and {\it divisorial hyperbolicity} notions, as well as a method for constructing new examples from existing ones.

\begin{The}\label{The:introd_hyperbolicity-examples} (i)\, The following two classes of compact complex manifolds $X$ consist exclusively of {\bf degenerate balanced} (hence also {\bf balanced hyperbolic}, hence also {\bf divisorially hyperbolic}) manifolds:

  \vspace{1ex}

  (a)\, the connected sums $X=\sharp_k(S^3\times S^3)$ of $k$ copies (with $k\geq 2$) of $S^3\times S^3$ endowed with the Friedman-Lu-Tian complex structure $J_k$ constructed via conifold transitions, where $S^3$ is the $3$-sphere;

  \vspace{1ex}

  (b)\, the Yachou manifolds $X=G/\Gamma$ arising as the quotient of any {\bf semi-simple} complex Lie group G by a lattice $\Gamma\subset G$.

  \vspace{2ex}

  (ii)\, If $X_1$ and $X_2$ are {\bf balanced hyperbolic} manifolds, so is $X_1\times X_2$.

  \vspace{2ex}

  (iii)\, The {\bf non-divisorially hyperbolic} manifolds include: all the complex projective spaces $\Proj^n$, all the complex tori $\C^n/\Gamma$, the $3$-dimensional Iwasawa manifold and all the $3$-dimensional Nakamura solvmanifolds.

\end{The}

A key feature of the Friedman-Lu-Tian manifolds $(X=\sharp_k(S^3\times S^3,\,J_k)$ and the Yachou manifolds $X=G/\Gamma$ is that their canonical bundle $K_X$ is {\it trivial}. By Kobayashi's Conjecture \ref{Conj:Kobayashi}, this is not expected to be the case in the classical context of Kobayashi/Brody hyperbolic compact complex manifolds $X$. We hope that examples of {\it projective} balanced hyperbolic or divisorially hyperbolic manifolds $X$ with $K_X$ {\it trivial} can be found in the future. Such an extension of the hyperbolicity theory into Calabi-Yau territory would be one of the spin-offs of our approach.

\vspace{2ex}

{\bf (II)} A fundamental problem in complex geometry is to prove positivity properties of various objects, notably the canonical bundle $K_X$, associated with a compact hyperbolic manifold $X$, as a way of emphasising the links between the complex analytic and the metric aspects of the theory.

In this vein, Gromov showed in [Gro91, Corollary 0.4.C] that $K_X$ is a {\it big} line bundle whenever $X$ is a compact K\"ahler hyperbolic manifold. Building on Gromov's result and on several classical results in birational geometry (including Mori's {\it Cone Theorem} implying that $K_X$ is nef whenever $X$ contains no rational curves, the Kawamata-Reid-Shokurov {\it Base-Point-Free Theorem} to the effect that $K_X$ is semi-ample whenever it is big and nef, and the {\it Relative Cone Theorem for log pairs}), Chen and Yang showed in [CY18, Theorem 2.11] that $K_X$ is even {\it ample} under the K\"ahler hyperbolicity assumption on the compact $X$. 

This answers affirmatively, in the special case of a compact K\"ahler hyperbolic manifold $X$, the following

\begin{Conj}(Kobayashi, see e.g. [CY18, Conjecture 2.8] or Lang's survey cited therein)\label{Conj:Kobayashi}

  If $X$ is a {\bf Kobayashi hyperbolic} compact complex manifold, its canonical bundle $K_X$ is {\bf ample}.

\end{Conj}

Our undertaking in $\S.$\ref{section:div-K_div-nef_classes} is motivated by a desire to prove positivity properties of $K_X$ under the weaker hyperbolicity assumptions on $X$ introduced in this paper, the {\it balanced hyperbolicity} and the {\it divisorial hyperbolicity}. Since there are quite a few {\it non-projective} and even {\it non-K\"ahler} compact complex manifolds $X$ that are hyperbolic in our two senses (see e.g. Theorem \ref{The:introd_hyperbolicity-examples}), $K_X$ cannot be positive in the usual big/ample senses for those $X$'s. Therefore, it seems natural to introduce positivity concepts relative to the complex codimension $1$ (rather than the usual complex dimension $1$) that hopefully match our codimension-$1$ hyperbolicity notions.

This is precisely what we propose in $\S.$\ref{section:div-K_div-nef_classes}. Given a compact complex $n$-dimensional manifold $X$, recall that the Bott-Chern and Aeppli cohomology groups of any bidegree $(p,\,q)$ of $X$ are classically defined, using the spaces $C^{r,\,s}(X) = C^{r,\,s}(X,\,\C)$ of smooth $\C$-valued $(r,\,s)$-forms on $X$, as \begin{eqnarray*}H^{p,\,q}_{BC}(X,\,\C) & = & \frac{\ker(\partial:C^{p,\,q}(X)\to C^{p+1,\,q}(X))\cap\ker(\bar\partial:C^{p,\,q}(X)\to C^{p,\,q+1}(X))}{\mbox{Im}\,(\partial\bar\partial:C^{p-1,\,q-1}(X)\to C^{p,\,q}(X))} \\
  H^{p,\,q}_A(X,\,\C) & = & \frac{\ker(\partial\bar\partial:C^{p,\,q}(X)\to C^{p+1,\,q+1}(X))}{\mbox{Im}\,(\partial:C^{p-1,\,q}(X)\to C^{p,\,q}(X)) + \mbox{Im}\,(\bar\partial:C^{p,\,q-1}(X)\to C^{p,\,q}(X))}.\end{eqnarray*} We will use the Serre-type duality (see e.g. [Sch07]): \begin{equation}\label{eqn:BC-A_duality}H^{1,\,1}_{BC}(X,\,\C)\times H^{n-1,\,n-1}_A(X,\,\C)\longrightarrow\C, \hspace{3ex} (\{u\}_{BC},\,\{v\}_A)\mapsto\{u\}_{BC}.\{v\}_A:=\int\limits_Xu\wedge v,\end{equation} as well as the {\it strongly Gauduchon (sG) cone} ${\cal SG}_X$ and the {\it Gauduchon cone} ${\cal G}_X$ of $X$ that were defined in [Pop15a] as: $${\cal SG}_X:=\bigg\{\{\omega^{n-1}\}_A\in H^{n-1,\,n-1}_A(X,\,\R)\,\mid\,\omega\hspace{1ex}\mbox{is an sG  metric}\hspace{1ex} \mbox{on}\hspace{1ex}X\bigg\}\subset H^{n-1,\,n-1}_A(X,\,\R);$$ $${\cal G}_X:=\bigg\{\{\omega^{n-1}\}_A\in H^{n-1,\,n-1}_A(X,\,\R)\,\mid\,\omega\hspace{1ex}\mbox{is a Gauduchon metric}\hspace{1ex} \mbox{on}\hspace{1ex}X\bigg\}\subset H^{n-1,\,n-1}_A(X,\,\R).$$ Recall that a Hermitian metric $\omega$ on $X$ is said to be a {\it Gauduchon metric} (cf. [Gau77a]), resp. a {\it strongly Gauduchon (sG) metric} (cf. [Pop13]), if $\partial\bar\partial\omega^{n-1}=0$, resp. if $\partial\omega^{n-1}\in\mbox{Im}\,\bar\partial$. Obviously, ${\cal SG}_X\subset{\cal G}_X$.

Now, given a real De Rham cohomology class $\{\alpha\}\in H^2_{DR}(X,\,\R)$ (not necessarily of type $(1,\,1)$), we say (see Definition \ref{Def:div-K-nef_classes}) that $\{\alpha\}$ is {\bf divisorially K\"ahler}, resp. {\bf divisorially nef}, if its image under the canonically defined map (see (\ref{eqn:P_n-1_map_def}) of Lemma \ref{Lem:P_map}): $$P:H^2_{DR}(X,\,\R)\longrightarrow H_A^{n-1,\,n-1}(X,\,\R), \hspace{3ex} \{\alpha\}_{DR}\longmapsto\{(\alpha^{n-1})^{n-1,\,n-1}\}_A,$$

\noindent lies in the {\it Gauduchon cone} ${\cal G}_X$ of $X$, respectively in the closure $\overline{\cal G}_X$ of this cone in $H^{n-1,\,n-1}_A(X,\,\R)$. 

We say that a $C^\infty$ complex line bundle $L$ on $X$ is {\it divisorially nef} if its first Chern class $c_1(L)$ is. An example of result in the special projective setting is the following immediate consequence of Propositions \ref{Prop:div-nef_classes_currents} and \ref{Prop:div-nef-cone_prop} (see $(4)$ of $\S.$\ref{subsection:div-nef_examp}): 

\begin{Prop}\label{Prop:introd_div-nef_char} Let $L$ be a holomorphic line bundle on an $n$-dimensional {\bf projective} manifold $X$. The following implication holds: $$L \hspace{2ex} \mbox{is {\bf divisorially nef}} \hspace{2ex} \implies \hspace{2ex} L^{n-1}.D\geq 0  \hspace{2ex} \mbox{for all effective divisors}  \hspace{2ex} D\geq 0 \hspace{2ex} \mbox{on} \hspace{1ex} X,$$ where $$L^{n-1}.D:=\int\limits_D\bigg(\frac{i}{2\pi}\Theta_h(L)\bigg)^{n-1}$$ and $(i/2\pi)\,\Theta_h(L)$ is the curvature form of $L$ with respect to any Hermitian fibre metric $h$.

\end{Prop}

If $L$ satisfies the last property above, we say that $L$ is {\it projectively divisorially nef}. This property is the divisorial analogue of the classical nefness property on projective manifolds $X$: $L$ is nef $\iff$ $L.C\geq 0$ for every curve $C\subset X$.

\vspace{2ex}

We also introduce the {\bf divisorially K\"ahler cone} ${\cal DK}_X$ and the {\bf divisorially nef cone} ${\cal DN}_X$ of $X$ in Definition \ref{Def:div-K-nef_classes} and discuss various properties of these notions in $\S.$\ref{subsection:div-nef_classes_proj-manif} and $\S.$\ref{subsection:div-nef_classes_arbitrary-manif}. In $\S.$\ref{subsection:div-nef_examp}, we point out examples of divisorially K\"ahler and divisorially nef cohomology classes.

Our hope is that we are able to take up the following problem in future work.

\begin{Question}\label{Question} Let $X$ be a compact complex manifold. If $X$ is {\bf balanced hyperbolic} or, more generally, {\bf divisorially hyperbolic}, does it follow that its canonical bundle $K_X$ is {\bf divisorially nef} or even {\bf divisorially K\"ahler}?
  
\end{Question}

\vspace{2ex}


\section{Balanced and divisorial hyperbolicity}\label{section:bal-div_hyperbolicity}

In this section, we introduce and discuss two hyperbolicity notions that generalise Gromov's K\"ahler hyperbolicity and the Kobayashi/Brody hyperbolicity respectively. 

\subsection{Balanced hyperbolic manifolds}\label{subsection:bal-hyp_manifolds}

Let $X$ be a compact complex manifold with $\mbox{dim}_\C X=n$. Fix an arbitrary Hermitian metric (i.e. a $C^\infty$ positive definite $(1,\,1)$-form) $\omega$ on $X$. Throughout the text, $\pi_X:\widetilde{X}\longrightarrow X$ will stand for the universal cover of $X$ and $\widetilde\omega=\pi_X^\star\omega$ will be the Hermitian metric on $\widetilde{X}$ that is the lift of $\omega$. Recall that a $C^\infty$ $k$-form $\alpha$ on $X$ is said to be $\widetilde{d}(\mbox{bounded})$ with respect to $\omega$ if $\pi_X^\star\alpha = d\beta$ on $\widetilde{X}$ for some $C^\infty$ $(k-1)$-form $\beta$ on $\widetilde{X}$ that is bounded w.r.t. $\widetilde\omega$. (See [Gro91].)

Recall two standard notions introduced by Gauduchon and Gromov respectively.

\vspace{1ex}

$(1)$\, The metric $\omega$ is said to be {\bf balanced} if $d\omega^{n-1}=0$. The manifold $X$ is said to be {\it balanced} if it carries a {\it balanced metric}. (See [Gau77b], where these metrics were called {\bf semi-K\"ahler}.)

\vspace{1ex}

$(2)$\, The metric $\omega$ is said to be {\bf K\"ahler hyperbolic} if $\omega$ is {\it K\"ahler} (i.e. $d\omega=0$) and $\widetilde{d}(\mbox{bounded})$ with respect to itself. The manifold $X$ is said to be {\it K\"ahler hyperbolic} if it carries a {\it K\"ahler hyperbolic metric}. (See [Gro91].) \\

\vspace{1ex}

The first notion that we introduce in this work combines the above two classical ones.

\begin{Def}\label{Def:bal-hyperbolic} Let $X$ be a compact complex manifold with $\mbox{dim}_\C X=n$. A Hermitian metric $\omega$ on $X$  is said to be {\bf balanced hyperbolic} if $\omega$ is balanced and $\omega^{n-1}$ is $\widetilde{d}(\mbox{bounded})$ with respect to $\omega$.

   The manifold $X$ is said to be balanced hyperbolic if it carries a {\it balanced hyperbolic metric}.

\end{Def}

Let us first notice the following implication: $$X \hspace{1ex} \mbox{is {\it K\"ahler hyperbolic}} \hspace{1ex} \implies \hspace{1ex} X \hspace{1ex} \mbox{is {\it balanced hyperbolic}}.$$

To see this, besides the obvious fact that every K\"ahler metric is balanced, we need the following

\begin{Lem}\label{Lem:powers_d-tilde-bounded} Let $(X,\,\omega)$ be a compact complex Hermitian manifold with $\mbox{dim}_\C X=n$. Let $k\in\{1,\dots , 2n\}$ and $\alpha\in C^\infty_k(X,\,\C)$. If $\alpha$ is $\widetilde{d}(\mbox{bounded})$ (with respect to $\omega$), then $\alpha^p$ is $\widetilde{d}(\mbox{bounded})$ (with respect to $\omega$) for every non-negative integer $p$.

\end{Lem}  

\noindent {\it Proof.} By the $\widetilde{d}(\mbox{boundedness})$ assumption on $\alpha$, $\pi_X^\star\alpha = d\beta$ on $\widetilde{X}$ for some smooth $\widetilde\omega$-bounded $(k-1)$-form $\beta$ on $\widetilde{X}$. Note that $d\beta$ is trivially $\widetilde\omega$-bounded on $\widetilde{X}$ since it equals $\pi_X^\star\alpha$ and $\alpha$ is $\omega$-bounded on $X$ thanks to $X$ being {\it compact}.

We get: $\pi_X^\star\alpha^p = d(\beta\wedge(d\beta)^{p-1})$ on $\widetilde{X}$, where both $\beta$ and $d\beta$ are $\widetilde\omega$-bounded, hence so is $\beta\wedge(d\beta)^{p-1}$.  \hfill $\qed$


\subsection{Divisorially hyperbolic manifolds}\label{subsection:div-hyp_manifolds}

We begin with a few preliminaries. Fix an arbitrary integer $n\geq 2$. For any $r>0$, let $B_r:=\{z\in\C^{n-1}\,\mid\,|z|<r\}$ and $S_r:=\{z\in\C^{n-1}\,\mid\,|z|=r\}$ stand for the open ball, resp. the sphere, of radius $r$ centred at $0\in\C^{n-1}$. Moreover, for any $(1,\,1)$-form $\gamma\geq 0$ on a complex manifold and any positive integer $p$, we will use the notation: $$\gamma_p:=\frac{\gamma^p}{p!}.$$

If $X$ is a compact complex manifold with $\mbox{dim}_\C X=n\geq 2$ and $\omega$ is a Hermitian metric on $X$, for any holomorphic map $f:\C^{n-1}\to X$ that is {\it non-degenerate} at some point $x\in\C^{n-1}$ (in the sense that its differential map $d_xf:\C^{n-1}\longrightarrow T_{f(x)}X$ at $x$ is of maximal rank), we consider the smooth $(1,\,1)$-form $f^\star\omega$ on $\C^{n-1}$. The assumptions made on $f$ imply that the differential map $d_zf$ is of maximal rank for every point $z\in\C^{n-1}\setminus\Sigma$, where $\Sigma\subset\C^{n-1}$ is an analytic subset. Thus, $f^\star\omega$ is $\geq 0$ on $\C^{n-1}$ and is $>0$ on $\C^{n-1}\setminus\Sigma$. Consequently, $f^\star\omega$ can be regarded as a {\it degenerate metric} on $\C^{n-1}$. Its degeneration locus, $\Sigma$, is empty if $f$ is non-degenerate at every point of $\C^{n-1}$, in which case $f^\star\omega$ is a genuine Hermitian metric on $\C^{n-1}$. However, in our case, $\Sigma$ will be non-empty in general, so $f^\star\omega$ will only be a genuine Hermitian metric on $\C^{n-1}\setminus\Sigma$.

For a holomorphic map $f:\C^{n-1}\to(X,\,\omega)$ in the above setting and for $r>0$, we consider the {\it $(\omega,\,f)$-volume} of the ball $B_r\subset\C^{n-1}$: $$\mbox{Vol}_{\omega,\,f}(B_r):=\int\limits_{B_r}f^\star\omega_{n-1}>0.$$

Meanwhile, for $z\in\C^{n-1}$, let $\tau(z):=|z|^2$ be its squared Euclidean norm. At every point $z\in\C^{n-1}\setminus\Sigma$, we have: \begin{eqnarray}\label{eqn:Hodge-star_area}\frac{d\tau}{|d\tau|_{f^\star\omega}}\wedge\star_{f^\star\omega}\bigg(\frac{d\tau}{|d\tau|_{f^\star\omega}}\bigg) = f^\star\omega_{n-1},\end{eqnarray} where $\star_{f^\star\omega}$ is the Hodge star operator induced by $f^\star\omega$. Thus, the $(2n-3)$-form $$d\sigma_{\omega,\,f}:=\star_{f^\star\omega}\bigg(\frac{d\tau}{|d\tau|_{f^\star\omega}}\bigg)$$ on $\C^{n-1}\setminus\Sigma$ is the area measure induced by $f^\star\omega$ on the spheres of $\C^{n-1}$. This means that its restriction \begin{eqnarray}\label{eqn:area-measure_t}d\sigma_{\omega,\,f,\,t}:=\bigg(\star_{f^\star\omega}\bigg(\frac{d\tau}{|d\tau|_{f^\star\omega}}\bigg)\bigg)_{|S_t}\end{eqnarray} is the area measure induced by the degenerate metric $f^\star\omega$ on the sphere $S_t=\{\tau(z)=t^2\}\subset\C^{n-1}$ for every $t>0$. In particular, the area of the sphere $S_r\subset\C^{n-1}$ w.r.t. $d\sigma_{\omega,\,f,\,r}$ is $$A_{\omega,\,f}(S_r)=\int\limits_{S_r}d\sigma_{\omega,\,f,\,r}>0,  \hspace{3ex} r>0.$$

\begin{Def}\label{Def:subexp-growth} Let $(X,\,\omega)$ be a {\bf compact} complex Hermitian manifold with $\mbox{dim}_\C X=n\geq 2$ and let $f:\C^{n-1}\to X$ be a holomorphic map that is {\bf non-degenerate} at some point $x\in\C^{n-1}$.

  We say that $f$ has {\bf subexponential growth} if the following two conditions are satisfied:

  \vspace{1ex}

  (i)\, there exist constants $C_1>0$ and $r_0> 0$ such that \begin{equation}\label{eqn:comparative-growth}\int\limits_{S_t}|d\tau|_{f^\star\omega}\,d\sigma_{\omega,\,f,\,t}\leq C_1t\,\mbox{Vol}_{\omega,\,f}(B_t),  \hspace{3ex} t>r_0;\end{equation}

\vspace{1ex} (ii)\, for every constant $C>0$, we have: \begin{equation}\label{eqn:subexp-growth}\limsup\limits_{b\to +\infty}\bigg(\frac{b}{C} - \log F(b)\bigg) = +\infty,\end{equation} where $$F(b):=\int\limits_0^b\mbox{Vol}_{\omega,\,f}(B_t)\,dt = \int\limits_0^b\bigg(\int\limits_{B_t}f^\star\omega_{n-1} \bigg)\,dt, \hspace{3ex} b>0.$$

\end{Def}

Note that (i), imposing a relative growth condition of the spheres $S_t$ w.r.t. the balls $B_t$ as measured by the degenerate metric $f^\star\omega$, is of a known type in this context. See, e.g. [dTh10]. The subexponential growth is expressed by condition (ii).

\vspace{2ex}

In a bid to shed light on the subexponential growth condition, we now spell out the very particular case where $f^\star\omega$ is the standard K\"ahler metric (i.e. the Euclidean metric) $\beta=(1/2)\,\sum\limits_{j=1}^{n-1} idz_j\wedge d\bar{z}_j$ of $\C^{n-1}$. It will come in handy when we discuss certain examples in $\S.$\ref{subsection:bal-hyp_examp}.

\begin{Lem}\label{Lem:subexp-growth_standard-metric} Let $$d\sigma_\beta:=\star_\beta\bigg(\frac{d\tau}{|d\tau|_\beta}\bigg)$$ be  the $(2n-3)$-form on $\C^{n-1}$ defining the area measure induced by $\beta$ on the spheres of $\C^{n-1}$. Then \begin{eqnarray}\label{eqn:area-measure_Euclidean}|d\tau|_\beta\,d\sigma_\beta = 2\sqrt\tau\,d\sigma_\beta  \hspace{5ex} \mbox{on} \hspace{2ex} \C^{n-1}.\end{eqnarray} In particular, \begin{equation}\label{eqn:integral_area-measure-tau_Euclidean}\int\limits_{S_t}|d\tau|_\beta\,d\sigma_\beta = 2A_{2n-3}\,t^{2n-2},  \hspace{3ex} t>0,\end{equation} where $A_{2n-3}$ is the area of the unit sphere $S_1\subset\C^{n-1}$ w.r.t. the measure $(d\sigma_\beta)_{|S_1}$ induced by the Euclidean metric $\beta$.

  In particular, any holomorphic map $f:\C^{n-1}\longrightarrow (X,\,\omega)$ such that $f^\star\omega = \beta$ has subexponential growth in the sense of Definition \ref{Def:subexp-growth}.

\end{Lem}

\noindent {\it Proof.} Since $d\tau = \partial\tau + \bar\partial\tau = \sum\limits_{j=1}^{n-1}\bar{z}_j\,dz_j + \sum\limits_{j=1}^{n-1}z_j\,d\bar{z}_j$ and $\langle dz_j,\,dz_k\rangle_\beta = \langle d\bar{z}_j,\,d\bar{z}_k\rangle_\beta = 2\delta_{jk}$, we get $|d\tau|_\beta^2 = 4|z|^2 = 4\tau$. This proves (\ref{eqn:area-measure_Euclidean}). Meanwhile, $\tau(z) = |z|^2 = t^2$ for $z\in S_t$, so we get (\ref{eqn:integral_area-measure-tau_Euclidean}).

On the other hand, $\mbox{Vol}_\beta(B_t) = V_{2n-2}\,t^{2n-2}$ for every $t>0$, so, when $f^\star\omega = \beta$, (\ref{eqn:comparative-growth}) amounts to $$2A_{2n-3}\,t^{2n-2} \leq C_1\,V_{2n-2}\,t^{2n-1},   \hspace{3ex} t>r_0,$$ which obviously holds for some constants $C_1, r_0>0$. Property (\ref{eqn:subexp-growth}) also holds in an obvious way. \hfill $\qed$

\vspace{2ex}

To further demystify condition (i) in Definition \ref{Def:subexp-growth}, we give an alternative expression for the integral on the sphere $S_t=\{|z|=t\}\subset\C^{n-1}$ featuring on the l.h.s. of (\ref{eqn:comparative-growth}) in terms of integrals on the ball $B_t=\{|z|<t\}\subset\C^{n-1}$. 

\begin{Lem}\label{Lem:alternative-formula_integral_growth} In the context of Definition \ref{Def:subexp-growth}, the following identities hold for all $t>0$: \begin{eqnarray}\label{eqn:alternative-formula_integral_growth}\int\limits_{S_t}|d\tau|_{f^\star\omega}\,d\sigma_{\omega,\,f,\,t} & = & 2\int\limits_{B_t}i\partial\bar\partial\tau\wedge f^\star\omega_{n-2} - \int\limits_{B_t}i(\bar\partial\tau - \partial\tau)\wedge d(f^\star\omega_{n-2}) \\
\nonumber    & = & 2\int\limits_{B_t}\Lambda_{f^\star\omega}(i\partial\bar\partial\tau)\,f^\star\omega_{n-1} - \int\limits_{B_t}i(\bar\partial\tau - \partial\tau)\wedge d(f^\star\omega_{n-2}),\end{eqnarray} where $\Lambda_{f^\star\omega}$ is the trace w.r.t. $f^\star\omega$ or, equivalently, the pointwise adjoint of the operator of multiplication by $f^\star\omega$, while $$i\partial\bar\partial\tau = i\partial\bar\partial|z|^2 = \sum\limits_{j=1}^{n-1}idz_j\wedge d\bar{z}_j:=\beta$$ is the standard metric of $\C^{n-1}$.

\end{Lem}

\noindent {\it Proof.} The pointwise identity $i\partial\bar\partial\tau\wedge(f^\star\omega)_{n-2} = \Lambda_{f^\star\omega}(i\partial\bar\partial\tau)\,(f^\star\omega)_{n-1}$ is standard on any $(n-1)$-dimensional complex manifold (which happens to be $\C^{n-1}$ in this case), so it suffices to prove the first equality in (\ref{eqn:alternative-formula_integral_growth}). 

We saw just above (\ref{eqn:area-measure_t}) that $|d\tau|_{f^\star\omega}\,d\sigma_{\omega,\,f} = \star_{f^\star\omega}(d\tau)$. Meanwhile, $d\tau = \partial\tau + \bar\partial\tau$ and the $1$-forms $\partial\tau$ and $\bar\partial\tau$ are primitive w.r.t. to any metric (in particular, w.r.t. $f^\star\omega$), as any $1$-form is. Consequently, the standard formula (cf. e.g. [Voi02, Proposition 6.29, p. 150]) for the Hodge star operator $\star = \star_\omega$ of any Hermitian metric $\omega$ applied to $\omega$-{\it primitive} forms $v$ of arbitrary bidegree $(p, \, q)$: \begin{eqnarray*}\star\, v = (-1)^{k(k+1)/2}\, i^{p-q}\,\omega_{n-p-q}\wedge v, \hspace{2ex} \mbox{where}\,\, k:=p+q,\end{eqnarray*} yields: $$\star_{f^\star\omega}(\partial\tau) = -i\partial\tau\wedge f^\star\omega_{n-2} \hspace{3ex} \mbox{and} \hspace{3ex} \star_{f^\star\omega}(\bar\partial\tau) = i\bar\partial\tau\wedge f^\star\omega_{n-2}.$$

Hence, we get the first equality below, where the second one follows from Stokes's theorem: \begin{eqnarray*}\label{eqn:alternative-formula_integral_growth_proof_1}\int\limits_{S_t}|d\tau|_{f^\star\omega}\,d\sigma_{\omega,\,f,\,t} & = & \int\limits_{S_t}i(\bar\partial\tau - \partial\tau)\wedge f^\star\omega_{n-2} = \int\limits_{B_t}d\bigg(i(\bar\partial\tau - \partial\tau)\wedge f^\star\omega_{n-2}\bigg) \\
 & = & \int\limits_{B_t}i\,d(\bar\partial\tau - \partial\tau)\wedge f^\star\omega_{n-2} - \int\limits_{B_t}i\,(\bar\partial\tau - \partial\tau)\wedge d(f^\star\omega_{n-2}),\end{eqnarray*} which is nothing but (\ref{eqn:alternative-formula_integral_growth}). \hfill $\qed$

\vspace{2ex}

Another immediate observation is that, due to $X$ being compact, we have

\begin{Lem}\label{Lem:subexp-growth_metric-independence} In the setting of Definition \ref{Def:subexp-growth}, the subexponential growth condition on $f$ is independent of the choice of Hermitian metric $\omega$ on $X$. 

\end{Lem}  

\noindent {\it Proof.} Let $\omega_1$ and $\omega_2$ be arbitrary Hermitian metrics on $X$. Since $X$ is compact, there exists a constant $A>0$ such that $(1/A)\,\omega_2\leq\omega_1\leq A\,\omega_2$ on $X$. Hence, $(1/A)\,f^\star\omega_2\leq f^\star\omega_1\leq A\,f^\star\omega_2$ on $\C^{n-1}$ for any holomorphic map $f:\C^{n-1}\to X$. The contention follows. \hfill $\qed$

\vspace{2ex}

Recall that a holomorphic map $f:\C^{n-1}\to(X,\,\omega)$ is standardly said to be {\it of finite order} if there exist constants $C_1,\,C_2,\,r_0>0$ such that \begin{equation}\label{eqn:finite-order}\mbox{Vol}_{\omega,\,f}(B_r)\leq C_1\,r^{C_2} \hspace{3ex} \mbox{for all} \hspace{1ex} r\geq r_0.\end{equation}

By the proof of Lemma \ref{Lem:subexp-growth_metric-independence}, $f$ being of finite order does not depend on the choice of Hermitian metric $\omega$ on $X$. Moreover, any $f$ of finite order satisfies condition (ii) in the definition \ref{Def:subexp-growth} of a subexponential growth. Furthermore, in the special case where $n-1=1$, it is a standard consequence of {\it Brody's Renormalisation Lemma} [Bro78, Lemma 2.1.] that any non-constant holomorphic map $f:\C\to X$ can be modified to a non-constant holomorphic map $\tilde{f}:\C\to X$ {\it of finite order}. (See e.g. [Lan87, Theorem 2.6., p. 72].)

However, one of the key differences between $\C$ and $\C^p$ with $p\geq 2$ is that a holomorphic map $\tilde{f}$ that is non-degenerate at some point (the higher dimensional analogue of the non-constancy of maps from $\C$) and has subexponential growth need not exist from $\C^{n-1}$ to a given compact $n$-dimensional $X$ when $n-1\geq 2$ even if a holomorphic map $f:\C^{n-1}\to X$ that is non-degenerate at some point exists. Based on this observation, we propose the following notion that generalises that of {\it Kobayashi/Brody hyperbolicity}.

\begin{Def}\label{Def:div-hyperbolic} Let $n\geq 2$ be an integer. An $n$-dimensional compact complex manifold $X$ is said to be {\bf divisorially hyperbolic} if there is no  holomorphic map $f:\C^{n-1}\longrightarrow X$ such that $f$ is {\bf non-degenerate} at some point $x\in\C^{n-1}$ and $f$ has {\bf subexponential growth} in the sense of Definition \ref{Def:subexp-growth}.

\end{Def} 


The following implication is obvious: $$X \hspace{1ex} \mbox{is {\it Kobayashi/Brody hyperbolic}} \hspace{1ex} \implies \hspace{1ex} X \hspace{1ex} \mbox{is {\it divisorially hyperbolic}}.$$ Indeed, if there is a holomorphic map $f:\C^{n-1}\longrightarrow X$ (of any growth) such that $f$ is non-degenerate at some point $x\in\C^{n-1}$, the restriction of $f$ to every complex line through $x$ in $\C^{n-1}$ is non-constant.

Meanwhile, the following implication is standard (see [Gro91, 0.3.B.]): $$X \hspace{1ex} \mbox{is {\it K\"ahler hyperbolic}} \hspace{1ex} \implies \hspace{1ex} X \hspace{1ex} \mbox{is {\it Kobayashi/Brody hyperbolic}}.$$

Taking our cue from the proof of Theorem 4.1 in [CY18], we now complete the diagram of implications in Theorem \ref{The:introd_hyperbolicity-implications} by proving its bottom row.

\begin{The}\label{The:bal-div-hyperbolic_implication} Every {\bf balanced hyperbolic} compact complex manifold is {\bf divisorially hyperbolic}.

\end{The}

\noindent {\it Proof.} Let $X$ be a compact complex manifold, with $\mbox{dim}_\C X=n$, equipped with a {\it balanced hyperbolic} metric $\omega$. This means that, if $\pi_X:\widetilde{X}\longrightarrow X$ is the universal cover of $X$, we have $$\pi_X^\star\omega^{n-1} = d\Gamma  \hspace{3ex} \mbox{on}\hspace{1ex} \widetilde{X},$$ where $\Gamma$ is an $\widetilde\omega$-bounded $C^\infty$ $(2n-3)$-form on $\widetilde{X}$ and $\widetilde\omega=\pi_X^\star\omega$ is the lift of the metric $\omega$ to $\widetilde{X}$.

Suppose there exists a holomorphic map $f:\C^{n-1}\longrightarrow X$ that is non-degenerate at some point $x\in\C^{n-1}$ and has subexponential growth in the sense of Definition \ref{Def:subexp-growth}. We will prove that $f^\star\omega^{n-1}=0$ on $\C^{n-1}$, in contradiction to the non-degeneracy assumption made on $f$ at $x$.

Since $\C^{n-1}$ is simply connected, there exists a lift $\widetilde{f}$ of $f$ to $\widetilde{X}$, namely a holomorphic map $\tilde{f}:\C^{n-1}\longrightarrow\widetilde{X}$ such that $f=\pi_X\circ\tilde{f}$. In particular, $d_x\tilde{f}$ is injective since $d_xf$ is.

The smooth $(n-1,\,n-1)$-form $f^\star\omega^{n-1}$ is $\geq 0$ on $\C^{n-1}$ and $>0$ on $\C^{n-1}\setminus\Sigma$, where $\Sigma\subset\C^{n-1}$ is the proper analytic subset of all points $z\in\C^{n-1}$ such that $d_zf$ is not of maximal rank. We have: $$f^\star\omega_{n-1} = \tilde{f}^\star(\pi_X^\star\omega^{n-1}) = d(\tilde{f}^\star\Gamma)  \hspace{3ex} \mbox{on}\hspace{1ex} \C^{n-1}.$$ 

 With respect to the degenerate metric $f^\star\omega$ on $\C^{n-1}$, we have the following

  \begin{Claim}\label{Claim:f-tilde-star_gamma_bounded} The $(2n-3)$-form $\tilde{f}^\star\Gamma$ is $(f^\star\omega)$-bounded on $\C^{n-1}$.

  \end{Claim}

  \noindent {\it Proof of Claim.} For any tangent vectors $v_1,\dots , v_{2n-3}$ in $\C^{n-1}$, we have: \begin{eqnarray*}|(\tilde{f}^\star\Gamma)(v_1,\dots , v_{2n-3})|^2 & = & |\Gamma(\tilde{f}_\star v_1,\dots , \tilde{f}_\star v_{2n-3})|^2 \stackrel{(a)}{\leq} C\,|\tilde{f}_\star v_1|^2_{\widetilde\omega}\dots |\tilde{f}_\star v_{2n-3}|^2_{\widetilde\omega} \\
    & = & C\,|v_1|^2_{\tilde{f}^\star\widetilde\omega}\dots |v_{2n-3}|^2_{\tilde{f}^\star\widetilde\omega} \stackrel{(b)}{=} C\,|v_1|^2_{f^\star\omega}\dots |v_{2n-3}|^2_{f^\star\omega},\end{eqnarray*} where $C>0$ is a constant independent of the $v_j$'s that exists such that inequality (a) holds thanks to the {\it $\widetilde\omega$-boundedness} of $\Gamma$ on $\widetilde{X}$, while (b) follows from $\tilde{f}^\star\widetilde\omega = f^\star\omega$.   \hfill $\qed$

  \vspace{2ex}

  \noindent {\it End of Proof of Theorem \ref{The:bal-div-hyperbolic_implication}.} We use the notation in the preliminaries of this $\S.$\ref{subsection:div-hyp_manifolds}.

\vspace{2ex}

 $\bullet$ On the one hand, we have $d\tau = 2t\,dt$ and \begin{eqnarray}\label{eqn:bal-div-hyperbolic_implication_proof_1_a}\mbox{Vol}_{\omega,\,f}(B_r) = \int\limits_{B_r}f^\star\omega_{n-1} = \int\limits_0^r\bigg(\int\limits_{S_t}d\mu_{\omega,\,f,\,t}\bigg)\,dt = \int_{B_r}d\mu_{\omega,\,f,\,t}\wedge\frac{d\tau}{2t},\end{eqnarray} where $d\mu_{\omega,\,f,\,t}$ is the positive measure on $S_t$ defined by $$\frac{1}{2t}\,d\mu_{\omega,\,f,\,t}\wedge (d\tau)_{|S_t} = (f^\star\omega_{n-1})_{|S_t},  \hspace{3ex} t>0.$$

Comparing with (\ref{eqn:Hodge-star_area}) and (\ref{eqn:area-measure_t}), this means that the measures $d\mu_{\omega,\,f,\,t}$ and $d\sigma_{\omega,\,f,\,t}$ on $S_t$ are related by \begin{eqnarray}\label{eqn:sigma-tau_relation}\frac{1}{2t}\,d\mu_{\omega,\,f,\,t} = \frac{1}{|d\tau|_{f^\star\omega}}\,d\sigma_{\omega,\,f,\,t}, \hspace{3ex} t>0.\end{eqnarray}

Now, the H\"older inequality yields: \begin{eqnarray*}\int_{S_t}\frac{1}{|d\tau|_{f^\star\omega}}\,d\sigma_{\omega,\,f,\,t}\geq\frac{A^2_{\omega,\,f}(S_t)}{\int_{S_t}|d\tau|_{f^\star\omega}\,d\sigma_{\omega,\,f,\,t}},\end{eqnarray*} so together with (\ref{eqn:bal-div-hyperbolic_implication_proof_1_a}) and (\ref{eqn:sigma-tau_relation}) this leads to: \begin{eqnarray}\label{eqn:bal-div-hyperbolic_implication_proof_1_b}\nonumber\mbox{Vol}_{\omega,\,f}(B_r) & = & \int\limits_0^r\bigg(\int\limits_{S_t}\frac{1}{2t}\,d\mu_{\omega,\,f,\,t}\bigg)\,d\tau = \int\limits_0^r\bigg(\int\limits_{S_t}\frac{1}{|d\tau|_{f^\star\omega}}\,d\sigma_{\omega,\,f,\,t}\bigg)\,d\tau \\
 & \geq & 2\,\int\limits_0^r\frac{A^2_{\omega,\,f}(S_t)}{\int_{S_t}|d\tau|_{f^\star\omega}\,d\sigma_{\omega,\,f,\,t}}\,t\,dt, \hspace{3ex} r>0.\end{eqnarray}

\vspace{2ex}

 $\bullet$ On the other hand, for every $r>0$, we have: \begin{eqnarray}\label{eqn:bal-div-hyperbolic_implication_proof_2}\mbox{Vol}_{\omega,\,f}(B_r) = \int\limits_{B_r} f^\star\omega_{n-1} = \int\limits_{B_r}d(\tilde{f}^\star\Gamma) = \int\limits_{S_r}\tilde{f}^\star\Gamma \stackrel{(a)}{\leq} C\,\int\limits_{S_r}d\sigma_{\omega,\,f} = C\,A_{\omega,\,f}(S_r),\end{eqnarray} where $C>0$ is a constant that exists such that inequality (a) holds thanks to Claim \ref{Claim:f-tilde-star_gamma_bounded}.

Putting (\ref{eqn:bal-div-hyperbolic_implication_proof_1_b}) and (\ref{eqn:bal-div-hyperbolic_implication_proof_2}) together, we get for every $r>r_0$: \begin{eqnarray}\label{eqn:bal-div-hyperbolic_implication_proof_3}\nonumber\mbox{Vol}_{\omega,\,f}(B_r) & \geq & \frac{2}{C^2}\,\int\limits_0^r\mbox{Vol}_{\omega,\,f}(B_t)\,\frac{t\,\mbox{Vol}_{\omega,\,f}(B_t)}{\int_{S_t}|d\tau|_{f^\star\omega}\,d\sigma_{\omega,\,f,\,t}} \,dt \\
   & \stackrel{(a)}{\geq} & \frac{2}{C_1\,C^2}\int\limits_{r_0}^r\mbox{Vol}_{\omega,\,f}(B_t)\,dt \stackrel{(b)}{:=} C_2\,F(r),\end{eqnarray} where (a) follows from the growth assumption (\ref{eqn:comparative-growth}) and (b) is the definition of a function $F:(r_0,\,+\infty)\longrightarrow(0,\,+\infty)$ with $C_2:=2/(C_1\,C^2)$.

  By taking the derivative of $F$, we get for every $r>r_0$: \begin{eqnarray*}F'(r) = \mbox{Vol}_{\omega,\,f}(B_r) \geq C_2\,F(r),\end{eqnarray*} where the last inequality is (\ref{eqn:bal-div-hyperbolic_implication_proof_3}). This amounts to \begin{eqnarray*}\frac{d}{dt}\bigg(\log F(t)\bigg) \geq C_2, \hspace{3ex} t>r_0.\end{eqnarray*} Integrating this over $t\in[a,\,b]$, with $r_0<a<b$ arbitrary, we get: \begin{eqnarray}\label{eqn:final_ineq_growth-cond}-\log F(a) \geq -\log F(b) + C_2\,(b-a), \hspace{3ex} r_0<a<b.\end{eqnarray}

  Now, fix an arbitrary $a>r_0$ and let $b\to +\infty$. Thanks to the {\it subexponential growth} assumption (\ref{eqn:subexp-growth}) made on $f$, there exists a sequence of reals $b_j\to +\infty$ such that the right-hand side of inequality (\ref{eqn:final_ineq_growth-cond}) for $b=b_j$ tends to $+\infty$ as $j\to +\infty$. This forces $F(a) = 0$ for every $a>r_0$, hence $\mbox{Vol}_{\omega,\,f}(B_r)=0$ for every $r>r_0$. This amounts to $f^\star\omega^{n-1}=0$ on $\C^{n-1}$, in contradiction to the non-degeneracy assumption made on $f$ at a point $x\in\C^{n-1}$.  \hfill $\qed$

\vspace{3ex}

   We now adapt in a straightforward way to our context the first part of the proof of [CY18, Proposition 2.11], where the non-existence of rational curves in compact K\"ahler hyperbolic manifolds was proved, and get the following analogous result.

  \begin{Prop}\label{Prop:bal-hyp_non-existence_Pn-1} Let $X$ be a compact complex manifold with $\mbox{dim}_\C X=n$. Suppose that $X$ carries a {\bf balanced hyperbolic} metric $\omega$. Then, there is no holomorphic map $f:\Proj^{n-1}\longrightarrow X$ such that $f$ is non-degenerate at some point $x\in\Proj^{n-1}$.

\end{Prop}

  \noindent {\it Proof.} Let $\pi_X^\star\omega_{n-1} = d\Gamma$ for some smooth $(2n-3)$-form $\Gamma$ on $\widetilde{X}$, where $\pi_X:\widetilde{X}\to X$ is the universal covering map of $X$. (We can even choose $\Gamma$ to be $\widetilde\omega$-bounded on $\widetilde{X}$, but we do not need this here.)

  Suppose there exists a holomorphic map $f:\Proj^{n-1}\longrightarrow X$ that is non-degenerate at some point. We will show that $f^\star\omega_{n-1}=0$ on $\Proj^{n-1}$, contradicting the non-degeneracy assumption on $f$.

 Let $\tilde{f}:\Proj^{n-1}\longrightarrow\widetilde{X}$ be a lift of $f$ to $\widetilde{X}$, namely a holomorphic map such that $f=\pi_X\circ\tilde{f}$. From $$f^\star\omega_{n-1} = \tilde{f}^\star(\pi_X^\star\omega_{n-1}) = d(\tilde{f}^\star\Gamma),$$ we get by integration: $$\int\limits_{\Proj^{n-1}}f^\star\omega_{n-1} = \int\limits_{\Proj^{n-1}}d(\tilde{f}^\star\Gamma) = 0,$$ where the last identity follows from Stokes's theorem.

Meanwhile, $f^\star\omega_{n-1}\geq 0$ at every point of $\Proj^{n-1}$. Therefore, $f^\star\omega_{n-1}=0$ on $\Proj^{n-1}$, a contradiction.  \hfill $\qed$

\subsection{Examples}\label{subsection:bal-hyp_examp} (I)\, The following definition was given in [Pop15a].

\begin{Def}\label{Def:deg-bal} Let $X$ be an $n$-dimensional complex manifold.

  A $C^\infty$ positive definite $(1,\,1)$-form $\omega$ on $X$ is said to be a {\bf degenerate balanced metric} if $\omega^{n-1}$ is $d$-exact. Any $X$ carrying such a metric is called a {\bf degenerate balanced manifold.}

\end{Def}

Degenerate balanced manifolds are characterised as follows.

\begin{Prop}\label{Prop:deg-bal_characterisation} ([Pop15a, Proposition 5.4]) Let $X$ be a compact complex manifold with $\mbox{dim}_C X=n$. The following statements are equivalent.

  \vspace{1ex}

  (i)\, The manifold $X$ is degenerate balanced.

\vspace{1ex}

  (ii)\, There exists no non-zero $d$-closed $(1, 1)$-current $T\geq 0$ on $X$.

\vspace{1ex}

  (iii)\, The Gauduchon cone of $X$ degenerates in the following sense: ${\cal G}_X = H_A^{n-1,\,n-1}(X,\,\R)$.

\end{Prop}

We are aware of two classes of {\bf degenerate balanced manifolds:}

\vspace{2ex}

(a)\, connected sums $X_k:=\sharp_k(S^3\times S^3)$ of $k\geq 2$ copies of $S^3\times S^3$, where $S^3$ is the unit sphere of $\R^4$ and each $X_k$ is endowed with the complex structure constructed by Friedman in [Fri89] and by Lu and Tian in [LT93] via {\it conifold transitions}. These complex structures were shown to be {\it balanced} in [FLY12]. Since $\mbox{dim}_\C X_k=3$ for every $k$ and since $H^4_{DR}(X_k,\,\C)=0$, any balanced metric $\omega_k$ on $X_k$ is necessarily {\it degenerate balanced}. Thus, $X_k$ is a {\it degenerate balanced manifold} for every $k\geq 2$. Note that every such $X_k$ is {\it simply connected}.

\vspace{1ex}

(b)\, quotients $X=G/\Gamma$ of a {\it semi-simple} complex Lie group $G$ by a lattice (i.e. a discrete co-compact subgroup) $\Gamma\subset G$. It was shown by Yachou in [Yac98, Propositions 17 and 18] that every left-invariant Hermitian metric on $G$ induces a {\it degenerate balanced} metric on $X$. Thus, any such $X$ (henceforth termed a {\bf Yachou manifold}) is a {\it degenerate balanced manifold}. Note that $X=G/\Gamma$ is {\it not simply connected} if $G$ is simply connected.

\vspace{2ex}

The immediate observation that provides the first class of examples of balanced hyperbolic manifolds is the following

\begin{Lem}\label{Lem:deg-bal_bal-hyp} Every {\bf degenerate balanced} compact complex manifold is {\bf balanced hyperbolic}.

\end{Lem}

\noindent {\it Proof.} If $\omega$ is a degenerate balanced metric on an $n$-dimensional $X$, then $\omega^{n-1} = d\Gamma$ for some smooth $(2n-3)$-form $\Gamma$ on $X$. Then, $\pi_X^\star(\omega^{n-1}) = d(\pi_X^\star\Gamma)$ on the universal covering manifold $\widetilde{X}$, while $\pi_X^\star\Gamma$ is $\widetilde\omega$-bounded on $\widetilde{X}$ since $\Gamma$ is $\omega$-bounded on the {\it compact} manifold $X$ as any smooth form is. As usual, $\pi_X:\widetilde{X}\longrightarrow X$ stands for the universal covering map and $\widetilde\omega:=\pi_X^\star\omega$.

Thus, $\omega^{n-1}$ is $\tilde{d}(\mbox{bounded})$ on $X$, so $\omega$ is a balanced hyperbolic metric.  \hfill $\qed$

\vspace{2ex}

Thus, for every compact complex manifold $X$, the following implications hold: \\

$X$ is {\it degenerate balanced} $\implies$ $X$ is {\it balanced hyperbolic} $\implies$  $X$ is {\it divisorially hyperbolic}.

\vspace{2ex}

Recall that a compact complex manifold $X$ is said to be {\it complex parallelisable} if its holomorphic tangent bundle $T^{1,\,0}X$ is trivial. By a result of Wang in [Wan54], $X$ is complex parallelisable if and only if $X$ is the compact quotient $X=G/\Gamma$ of a simply connected, connected {\bf complex} Lie group $G$ by a discrete subgroup $\Gamma\subset G$.

Meanwhile, it is standard that no {\it complex} Lie group $G$ is Brody hyperbolic. Indeed, the complex one-parameter subgroup generated by any given element $\xi$ in the Lie algebra of $G$ provides an example of an entire curve in $G$. (In other words, take any tangent vector $\xi$ of type $(1,\,0)$ in the tangent space $T_eG$ at the identity element $e\in G$ seen as the Lie algebra $\fg$ of $G$, then compose the linear map from $\C$ to $\fg$ that maps $1$ to $\xi$ with the exponential map $\exp:\fg\to G$, which is holomorphic since $G$ is a complex Lie group, to get an entire curve in $G$.) Together with Wang's result mentioned above, this shows that no {\it complex parallelisable} compact complex manifold $X$ is {\it Brody hyperbolic} (or, equivalently, since $X$ is compact, {\it Kobayashi hyperbolic}).

Now, note that the Yachou manifolds $X=G/\Gamma$ mentioned above are {\it complex parallelisable} manifolds since $G$ is a {\it complex} Lie group. So, they are {\it not Kobayashi hyperbolic}. In particular, they are {\it not K\"ahler hyperbolic}. However, they are {\it degenerate balanced}, hence also {\it balanced hyperbolic} (by Lemma \ref{Lem:deg-bal_bal-hyp}), hence also {\it divisorially hyperbolic} (by Theorem \ref{The:bal-div-hyperbolic_implication}). On the other hand, the Yachou manifolds $X=G/\Gamma$ are not K\"ahler (since, for example, they do not even support non-zero $d$-closed positive $(1,\,1)$-currents, by Proposition \ref{Prop:deg-bal_characterisation}). Hence, we get the following observation showing that the notions of {\it balanced hyperbolic} manifolds and {\it divisorially hyperbolic} manifolds are new and propose a hyperbolicity theory in the possibly non-K\"ahler context.

\begin{Prop}\label{Prop:div-hyperbolic_not-Kob-hyperbolic} There exist compact complex non-K\"ahler manifolds that are {\bf balanced hyperbolic} but are not {\bf Kobayashi hyperbolic}.

\end{Prop}  

\vspace{2ex}

It seems natural to ask the following

\begin{Question}\label{Question:compact-quotients_div-hyp} Which compact quotients $X=G/\Gamma$ of a {\bf complex} Lie group $G$ by a lattice $\Gamma$ are {\bf balanced hyperbolic} or, at least, {\bf divisorially hyperbolic}?

\end{Question}  

We know from [Yac98] that all these quotients are even {\it degenerate balanced} (hence also {\it balanced hyperbolic}, hence also {\it divisorially hyperbolic}) when $G$ is {\it semi-simple}. On the other hand, there is always a {\it holomorphic map} $f:\C^{n-1}\to X=G/\Gamma$, non-degenerate at some point $x\in\C^{n-1}$, whenever $G$ is an $n$-dimensional complex Lie group and $\Gamma\subset G$ is a discrete co-compact subgroup. Indeed, let $\xi_1,\dots , \xi_{n-1}\in T_eG=\fg$ be $\C$-linearly independent vectors of type $(1,\,0)$ in the Lie algebra of $G$ and let $h:\C^{n-1}\to\fg$ be the $\C$-linear map that takes the vectors $e_1,\dots , e_{n-1}$ forming the canonical basis of $\C^{n-1}$ to $\xi_1,\dots , \xi_{n-1}$ respectively. The desired map $f:\C^{n-1}\to X=G/\Gamma$ is obtained by composing $h$ with the exponential map $\exp:\fg\to G$ (which is holomorphic, due to $G$ being a complex Lie group, and non-degenerate at least at $0\in\fg$, hence at least on a neighbourhood of it, since its differential map at $0$ is the identity map) and with the projection map $G\to G/\Gamma$.

Thus, part of Question \ref{Question:compact-quotients_div-hyp} reduces to determining the complex Lie groups $G$ and their lattices $\Gamma$ for which no map as above that also has {\it subexponential growth} in the sense of Definition \ref{Def:subexp-growth} exists. Meanwhile, we point to (a), (b), (c) under (VI) in this $\S.$\ref{subsection:bal-hyp_examp} for examples of non-hyperbolic compact quotients $G/\Gamma$ of a complex Lie group by a lattice.

\vspace{2ex}

(II)\, The other obvious class of {\it balanced hyperbolic manifolds} consists of all the {\it K\"ahler hyperbolic manifolds}. (See $\S.$\ref{section:bal-div_hyperbolicity}.)

\vspace{2ex}

(III)\, We shall now point out examples of {\it balanced hyperbolic} manifolds that are neither degenerate balanced, nor K\"ahler hyperbolic. We first recall the following result of Michelsohn.

\begin{Prop}\label{Prop:product_balanced} ([Mic83, Proposition 1.9]) Let $X$ and $Y$ be complex manifolds.

\vspace{1ex}

(i)\, If $X$ and $Y$ are balanced, the product manifold $X\times Y$ is balanced.

\vspace{1ex}

(ii)\, Let $\sigma_X$ and $\sigma_Y$ be the projections of $X\times Y$ onto $X$, resp. $Y$. If $\omega_X$ and $\omega_Y$ are balanced metrics on $X$, resp. $Y$, the induced product metric $\omega=\sigma_X^\star\omega_X + \sigma_Y^\star\omega_Y$ is a balanced metric on $X\times Y$.

\end{Prop}

\vspace{2ex}

Using this, we notice the following simple way of producing new balanced hyperbolic manifolds from existing ones. 

\begin{Prop}\label{Prop:product_bal-hyperbolic} The Cartesian product of {\bf balanced hyperbolic} manifolds is {\bf balanced hyperbolic}.

\end{Prop}

 \noindent {\it Proof.} Let $(X_1,\omega_1)$ and $(X_2,\omega_2)$ be balanced hyperbolic manifolds of respective dimensions $n$ and $m$, and let $\pi_1:\widetilde{X_1}\longrightarrow X_1$, $\pi_2:\widetilde{X_2}\longrightarrow X_2$ be their respective universal covers. By hypothesis, we have:

 \vspace{1ex}
 
 $\bullet$ $\omega_1^{n-1}$ is $\tilde{d}$\emph{(bounded)} on $(X_1,\omega_1)$, so there exists an $\widetilde\omega_1$-bounded $(2n-3)$-form $\Theta_1$ on $\widetilde{X}_1$ such that $\pi_1^*(\omega_1^{n-1})=d\Theta_1$, where $\widetilde\omega_1:=\pi_1^\star\omega_1$;

\vspace{1ex}
 
  $\bullet$ $\omega_2^{m-1}$ is $\tilde{d}$\emph{(bounded)} on $(X_2,\omega_2)$, so there exists an $\widetilde\omega_2$-bounded $(2m-3)$-form $\Theta_2$ on $\widetilde{X}_2$ such that $\pi_2^*(\omega_2^{m-1})=d\Theta_2$, where $\widetilde\omega_2:=\pi_2^\star\omega_2$.
     
The product map $$\pi:=\pi_1\times\pi_2:\widetilde{X_1}\times\widetilde{X_2}\longrightarrow X_1\times X_2 $$ is the universal cover of $X_1\times X_2$. Meanwhile, by (ii) of Proposition \ref{Prop:product_balanced}, the product metric $$\omega=\sigma_1^\star\omega_1+\sigma_2^\star\omega_2$$ on $X_1\times X_2$ is  balanced, where $\sigma_1:X_1\times X_2\to X_1$ and $\sigma_2:X_1\times X_2\to X_2$ are the projections on the two factors. From the equality $$\omega^{n+m-1}={n+m-1 \choose n-1}\,\sigma_1^\star\omega_1^{n-1}\wedge\sigma_2^\star\omega_2^m+{n+m-1 \choose n}\,\sigma_1^\star\omega_1^{n}\wedge\sigma_2^\star\omega_2^{m-1}$$ on $X_1\times X_2$, we infer the following equalities on $\widetilde{X_1}\times \widetilde{X_2}$:

\begin{eqnarray}\label{eqn:product_bal-hyperbolic_proof_1}\nonumber\pi^\star(\omega^{n+m-1})&=&{n+m-1 \choose n-1}\,\pi^\star(\sigma_1^\star(\omega_1^{n-1}))\wedge\pi^\star(\sigma_2^\star(\omega_2^m))+{n+m-1 \choose n}\,\pi^\star(\sigma_1^\star(\omega_1^n))\wedge\pi^\star(\sigma_2^\star(\omega_2^{m-1}))\\
\nonumber &\stackrel{(a)}{=}&{n+m-1 \choose n-1}\,\widetilde{\sigma_1}^\star(\pi_1^\star(\omega_1^{n-1}))\wedge\widetilde{\sigma_2}^\star(\pi_2^\star(\omega_2^m))+{n+m-1 \choose n}\,\widetilde{\sigma_1}^\star(\pi_1^\star(\omega_1^n))\wedge\widetilde{\sigma_2}^\star(\pi_2^\star(\omega_2^{m-1}))\\
\nonumber &=&{n+m-1 \choose n-1}\,\widetilde{\sigma_1}^\star(d\Theta_1)\wedge\widetilde{\sigma_2}^\star(\pi_2^\star(\omega_2^m))+{n+m-1 \choose n}\,\widetilde{\sigma_1}^\star(\pi_1^\star(\omega_1^n))\wedge\widetilde{\sigma_2}^\star(d\Theta_2)\\
\nonumber &=&{n+m-1 \choose n-1}\,d(\widetilde{\sigma_1}^\star \Theta_1\wedge\widetilde{\sigma_2}^\star(\pi_2^\star(\omega_2^m)))+{n+m-1 \choose n}\,d\left(\widetilde{\sigma_1}^\star(\pi_1^\star(\omega_1^n))\wedge\widetilde{\sigma_2}^\star\Theta_2\right)\\
&=&d\left[{n+m-1 \choose n-1}\,\widetilde{\sigma_1}^\star\Theta_1\wedge\widetilde{\sigma_2}^\star(\pi_2^\star(\omega_2^m))+{n+m-1 \choose n}\,\widetilde{\sigma_1}^\star(\pi_1^\star(\omega_1^n))\wedge\widetilde{\sigma_2}^\star \Theta_2\right]
\end{eqnarray} where $\widetilde{\sigma_1}:\widetilde{X_1}\times \widetilde{X_2}\to \widetilde{X_1}$ and $\widetilde{\sigma}_2:\widetilde{X_1}\times \widetilde{X_2}\to \widetilde{X_2}$ are the projections on the two factors. Note that the equalities $\sigma_j\circ\pi = \pi_j\circ\widetilde{\sigma_j}$ for $j=1,\,2$ were used to get equality (a) in (\ref{eqn:product_bal-hyperbolic_proof_1}).

Now, for every $j\in\{1,\,2\}$, $\Theta_j$ is $\widetilde\omega_j$-bounded on $\widetilde{X}_j$. Therefore, $\widetilde{\sigma_j}^\star\Theta_j$ is $\widetilde{\sigma_j}^\star\widetilde\omega_j$-bounded, hence also $\widetilde\omega$-bounded, on $\widetilde{X_1}\times \widetilde{X_2}$, where $\widetilde\omega$ is the product metric $$\widetilde\omega=\widetilde\sigma_1^\star\widetilde\omega_1+\widetilde\sigma_2^\star\widetilde\omega_2$$ on $\widetilde{X}_1\times\widetilde{X}_2$. We infer that the forms $\widetilde{\sigma_1}^\star \Theta_1\wedge\widetilde{\sigma_2}^\star(\pi_2^\star(\omega_2^m))$ and $\widetilde{\sigma_1}^\star(\pi_1^\star(\omega_1^n))\wedge\widetilde{\sigma_2}^\star \Theta_2$ are $\widetilde\omega$-bounded on $\widetilde{X_1}\times \widetilde{X_2}$, hence so is their linear combination featuring in the $d$-potential of the form on the r.h.s. of the last line in (\ref{eqn:product_bal-hyperbolic_proof_1}).

 Consequently, the form $\omega^{n+m-1}$ is $\tilde{d}$\emph{(bounded)} on $(X_1\times X_2,\,\omega)$. (Note that $\pi^\star\omega= \widetilde\omega$.) This means that the metric $\omega$ of the $(n+m)$-dimensional complex manifold $X_1\times X_2$ is balanced hyperbolic. Hence, the manifold $X_1\times X_2$ is balanced hyperbolic. \hfill $\qed$\\

 In particular, using the above result, we can construct examples of {\it non-K\"ahler} balanced hyperbolic manifolds.

\begin{Cor}\label{Cor:product_bal-hyperbolic} If $X_1$ is any {\bf degenerate balanced} manifold and $X_2$ is any {\bf K\"ahler hyperbolic} manifold, $X_1\times X_2$ is a {\bf balanced hyperbolic} manifold that need not be either degenerate balanced, or K\"ahler hyperbolic, or even K\"ahler.

\end{Cor}

On the other hand, we can also construct examples of {\it K\"ahler} balanced hyperbolic manifolds that are neither K\"ahler hyperbolic, nor degenerate balanced.

\begin{Prop}\label{Prop:Kaehler_product_bal-hyperbolic}
  Let $X_1$ be a {\bf K\"ahler hyperbolic} manifold with $\mbox{dim}_\C X_1 =n>1$ and let $X_2$ be a {\bf compact K\"ahler} manifold. Then, $X_1\times X_2$ is a {\bf balanced hyperbolic} K\"ahler manifold.
\end{Prop}

\noindent {\it Proof.} Let $\omega_1$ be a K\"ahler hyperbolic metric on $X_1$, $\omega_2$ a K\"ahler metric on $X_2$ and $m=\mbox{dim}_\C X_2$. We will keep the notation used in the proof of Proposition \ref{Prop:product_bal-hyperbolic}, except for the differences that will be pointed out.

Since $\omega_1$ is $\tilde{d}$(\emph{bounded}), so are $\omega_1^{n-1}$  and $\omega_1^n$, by Lemma \ref{Lem:powers_d-tilde-bounded}. Thus, there exist $\widetilde\omega_1$-bounded forms $\Theta_1$ and $\Gamma_1$ on $\widetilde{X}_1$, of respective degrees $(2n-3)$ and $(2n-1)$, such that $$\pi_1^*(\omega_1^{n-1})=d\Theta_1 \hspace{3ex}\mbox{and}\hspace{3ex} \pi_1^*(\omega_1^{n})=d\Gamma_1.$$

The only differences from the proof of Proposition \ref{Prop:product_bal-hyperbolic} are the disappearance of $\Theta_2$ and the appearance of $\Gamma_1$, together with the different properties that $\omega_1$ and $\omega_2$ now have. Running the equalities (\ref{eqn:product_bal-hyperbolic_proof_1}) with these differences incorporated, we get on $\widetilde{X_1}\times \widetilde{X_2}$: \begin{eqnarray*}\pi^\star(\omega^{n+m-1})&=& d \left[{n+m-1 \choose n-1}\,\widetilde{\sigma_1}^\star \Theta_1\wedge\widetilde{\sigma_2}^\star(\pi_2^\star(\omega_2^m))+{n+m-1 \choose n}\,\widetilde{\sigma_1}^\star\Gamma_1\wedge\widetilde{\sigma_2}^\star(\pi_2^\star(\omega_2^{m-1}))\right]\end{eqnarray*} after using the fact that $d\omega_2=0$ (the K\"ahler assumption on $\omega_2$). 

We conclude in the same way as in the proof of Proposition \ref{Prop:product_bal-hyperbolic} that the $d$-potential on the r.h.s. of the last line above is $\widetilde\omega$-bounded on $\widetilde{X_1}\times \widetilde{X_2}$. Thus,  the form $\omega^{n+m-1}$ is $\tilde{d}$\emph{(bounded)} on $(X_1\times X_2,\,\omega)$, so $\omega$ is a balanced hyperbolic metric on $X_1\times X_2$. \hfill $\qed$

 \vspace{2ex}

 (IV)\, We now discuss in some detail the case of the semi-simple complex Lie group $G=SL(2,\,\C)$, where several of the above constructions can be described explicitly.

 As a complex manifold, $G=SL(2,\,\C)$ is of complex dimension $3$. Its complex structure is described by three holomorphic $(1,\,0)$-forms $\alpha,\beta,\gamma$ that satisfy the structure equations: \begin{equation}\label{eqn:structure-eq_SL2C}d\alpha = \beta\wedge\gamma, \hspace{3ex} d\beta = \gamma\wedge\alpha, \hspace{3ex} d\gamma = \alpha\wedge\beta.\end{equation} Moreover, the dual of the Lie algebra $\fg = T_eG$ of $G$ is generated, as an $\R$-vector space, by these forms and their conjugates: $$(T_eG)^\star = \langle\alpha,\,\beta,\,\gamma,\,\overline\alpha,\,\overline\beta,\,\overline\gamma\rangle.$$

 The $C^\infty$ positive definite $(1,\,1)$-form $$\omega:=\frac{i}{2}\alpha\wedge\overline\alpha + \frac{i}{2}\beta\wedge\overline\beta + \frac{i}{2}\gamma\wedge\overline\gamma$$ defines a left-invariant (under the action of $G$ on itself) Hermitian metric on $G$. From this, using (\ref{eqn:structure-eq_SL2C}), we get $$\omega^2 = \frac{1}{2}\,d(\alpha\wedge d\overline\alpha + \beta\wedge d\overline\beta + \gamma\wedge d\overline\gamma)\in\mbox{Im}\,d.$$ So, $\omega$ is a {\it degenerate balanced} metric on $G$ (see Definition \ref{Def:deg-bal}). Since it is left-invariant under the $G$-action, $\omega$ descends to a {\it degenerate balanced metric} on the compact quotient $X=G/\Gamma$ of $G$ by any lattice $\Gamma$. In particular, this example illustrates Yachou's result [Yac98, Propositions 17 and 18] in the special case of $G=SL(2,\,\C)$.

 Now, consider the holomorphic map \begin{equation}\label{eqn:hol-map_C2-SL2C}f:\C^2\to G=SL(2,\,\C), \hspace{3ex} f(z_1,\,z_2) = \begin{pmatrix}e^{z_1} & z_2 \\
     0 & e^{-z_1}\end{pmatrix}.\end{equation} This map is non-degenerate at every point $z=(z_1,\,z_2)\in\C^2$, as can be seen at once. However, $f$ is {\it not of subexponential growth} in the sense of Definition \ref{Def:subexp-growth}, as we will now see. Actually, there is no non-degenerate holomorphic map $g:\C^2\to X=G/\Gamma$ of subexponential growth thanks to $X$ being {\it degenerate balanced} (hence also {\it balanced hyperbolic}) and to our Theorem \ref{The:bal-div-hyperbolic_implication}.

 \begin{Lem}\footnote{The authors are very grateful to L. Ugarte for pointing out to them the calculations leading to the result of this lemma and most of the discussion of $SL(2,\,\C)$ preceding it under (IV).}\label{Lem:SL2C_f-star-omega_computation} In $\C^2$, we have \begin{eqnarray*}f^\star\omega & = & (|z_2|^2 e^{-2\mbox{Re}\,z_1} + 2)\,idz_1\wedge d\bar{z}_1 + e^{-2\mbox{Re}\,z_1}\,idz_2\wedge d\bar{z}_2 \\
  & + & z_2\,e^{-2\mbox{Re}\,z_1}\,idz_1\wedge d\bar{z}_2 + \bar{z}_2\,e^{-2\mbox{Re}\,z_1}\,idz_2\wedge d\bar{z}_1.\end{eqnarray*}

 \end{Lem}

 \noindent {\it Proof.} $\bullet$ {\it Calculations at $(0,\,0)\in\C^2$}. \begin{eqnarray*}f_\star\bigg(\frac{\partial}{\partial z_1}_{|(0,\,0)}\bigg) & = & \frac{d}{dt}_{|t=0} f(t,\,0) = \frac{d}{dt}_{|t=0}\begin{pmatrix}e^t & 0 \\
     0 & e^{-t}\end{pmatrix} = \begin{pmatrix}1 & 0 \\
     0 & -1\end{pmatrix}:=H \\
     f_\star\bigg(\frac{\partial}{\partial z_2}_{|(0,\,0)}\bigg) & = & \frac{d}{dt}_{|t=0} f(0,\,t) = \frac{d}{dt}_{|t=0}\begin{pmatrix}1 & t \\
     0 & 1\end{pmatrix} = \begin{pmatrix}0 & 1 \\
     0 & 0\end{pmatrix}:=X.\end{eqnarray*} Note that $H,X\in T_{I_2}SL(2,\,\C) = sl(2,\,\C).$ A basis of the Lie algebra $sl(2,\,\C)$ is given by $\{H, X, Y\}$, where $$Y:= \begin{pmatrix}0 & 0 \\
   1 & 0\end{pmatrix}.$$ The Lie brackets linking the elements of this basis are \begin{eqnarray}\label{eqn:Lie-brackets_basis_1}[X,\,Y] = H; \hspace{3ex} [Y,\,H] = 2Y; \hspace{3ex}  [X,\,H] = -2X.\end{eqnarray}

   Rather than being dual to the basis $\{H, X, Y\}$, the basis $\{\alpha,\,\beta,\,\gamma\}$ of left-invariant $(1,\,0)$-forms that satisfy equations (\ref{eqn:structure-eq_SL2C}) is dual to the following basis of tangent vectors of type $(1,\,0)$ at $I_2\in SL(2,\,\C)$: \begin{eqnarray*}A:=\frac{i}{2}\,(X+Y) = \frac{i}{2}\,\begin{pmatrix}0 & 1 \\
     1 & 0\end{pmatrix}, \hspace{2ex} B:=\frac{1}{2}\,(X-Y) = \frac{1}{2}\,\begin{pmatrix}0 & 1 \\
     -1 & 0\end{pmatrix}, \hspace{2ex} C:=\frac{i}{2}\,H = \frac{i}{2}\,\begin{pmatrix}1 & 0 \\
     0 & -1\end{pmatrix}.\end{eqnarray*} This amounts to \begin{eqnarray*}X = -iA+B, \hspace{3ex} Y = -iA-B, \hspace{3ex} H=-2iC.\end{eqnarray*} From (\ref{eqn:Lie-brackets_basis_1}), we get \begin{eqnarray}\label{eqn:Lie-brackets_basis_2}[A,\,B] = -C; \hspace{3ex} [A,\,C] = B; \hspace{3ex}  [B,\,C] = -A.\end{eqnarray}
 (To see this, observe, for example, that $[A,\,B] = -(i/2)\,[X,\,Y] = -(i/2)\,H = -C$.)

\vspace{2ex}

$\bullet$ {\it Calculations at an arbitrary point $(z_1^0,\,z_2^0)\in\C^2$}.   

\vspace{1ex}

Since $T_{I_2}SL(2,\,\C) = \langle H,\,X,\,Y\rangle = \langle A,\,B,\,C\rangle$, for every $g\in SL(2,\,\C)$, the tangent space at $g$ is generated as $$T_gSL(2,\,\C) = \langle (L_g)_\star H,\,(L_g)_\star X,\,(L_g)_\star Y\rangle = \langle (L_g)_\star A,\,(L_g)_\star B,\,(L_g)_\star C\rangle,$$ where $L_g:G\to G$ is the left translation by $g$, namely $L_g(h) = gh$ for every $h\in G$. Since $L_g$ is a linear map, $(L_g)_\star = L_g$. Therefore, for every $$g=\begin{pmatrix}a & b \\
     c & d\end{pmatrix}\in SL(2,\,\C) \hspace{3ex} (\mbox{with}\hspace{1ex} ad-bc=1),$$ we get:  \begin{eqnarray*}(L_g)_\star A & = & (L_g) A = \frac{i}{2}\,\begin{pmatrix}a & b \\
     c & d\end{pmatrix} \begin{pmatrix}0 & 1 \\
     1 & 0\end{pmatrix} = \frac{i}{2}\,\begin{pmatrix}b & a \\
       d & c\end{pmatrix}:= A_g \\
     (L_g)_\star B & = & (L_g) B = \frac{1}{2}\,\begin{pmatrix}a & b \\
     c & d\end{pmatrix} \begin{pmatrix}0 & 1 \\
     -1 & 0\end{pmatrix} = \frac{1}{2}\,\begin{pmatrix}-b & a \\
       -d & c\end{pmatrix}:= B_g \\
     (L_g)_\star C & = & (L_g) C = \frac{i}{2}\,\begin{pmatrix}a & b \\
     c & d\end{pmatrix} \begin{pmatrix}1 & 0 \\
     0 & -1\end{pmatrix} = \frac{i}{2}\,\begin{pmatrix}a & -b \\
       c & -d\end{pmatrix}:= C_g.\end{eqnarray*}

 We now fix an arbitrary point $(z_1^0,\,z_2^0)\in\C^2$ and we let  \begin{eqnarray*} g_0:=f(z_1^0,\,z_2^0) = \begin{pmatrix}e^{z_1^0} & z_2^0 \\
     0 & e^{-z_1^0}\end{pmatrix}\in SL(2,\,\C).\end{eqnarray*} (Thus, for $g_0$, we have: $a=e^{z_1^0}$, $b=z_2^0$, $c=0$, $d=e^{-z_1^0}$.) We get:

\begin{eqnarray}\label{eqn:f_star_z_1}\nonumber f_\star\bigg(\frac{\partial}{\partial z_1}_{|(z_1^0,\,z_2^0)}\bigg) & = & \frac{d}{dt}_{|t=0} f((z_1^0,\,z_2^0)+ (t,\,0)) = \frac{d}{dt}_{|t=0}\begin{pmatrix}e^{z_1^0 + t} & z_2^0 \\
    0 & e^{-z_1^0-t}\end{pmatrix} \\
    \nonumber & = & \begin{pmatrix}e^{z_1^0} & 0 \\
    0 & -e^{-z_1^0}\end{pmatrix} = \begin{pmatrix}e^{z_1^0} & -z_2^0 \\
    0 & -e^{-z_1^0}\end{pmatrix} + e^{-z_1^0}z_2^0\begin{pmatrix}0 & e^{z_1^0} \\
  0 & 0\end{pmatrix} \\
   & = & -2iC_{g_0} + e^{-z_1^0}z_2^0\,(-iA_{g_0} + B_{g_0}).\end{eqnarray} Similarly, we get \begin{eqnarray}\label{eqn:f_star_z_2}\nonumber f_\star\bigg(\frac{\partial}{\partial z_2}_{|(z_1^0,\,z_2^0)}\bigg) & = & \frac{d}{dt}_{|t=0} f((z_1^0,\,z_2^0)+ (0,\,t)) = \frac{d}{dt}_{|t=0}\begin{pmatrix}e^{z_1^0} & z_2^0 + t \\
    0 & e^{-z_1^0}\end{pmatrix} \\
     & = & \begin{pmatrix}0 & 1 \\
    0 & 0\end{pmatrix} = e^{-z_1^0}\,\begin{pmatrix}0 & e^{z_1^0} \\
    0 & 0\end{pmatrix} = e^{-z_1^0}\,(-iA_{g_0} + B_{g_0}).\end{eqnarray}

We now use the general formula $(f^\star\omega)(V,\,W) = \omega(f_\star V,\,f_\star W)$ (for all vector fields $V,W$) to deduce expressions for $(f^\star\omega)(\partial/\partial z_j,\,\partial/\partial\bar{z}_k)$ (for all $j,k=1,2$) from (\ref{eqn:f_star_z_1}) and (\ref{eqn:f_star_z_2}).

From (\ref{eqn:f_star_z_1}), we get: \begin{eqnarray}\nonumber\label{eqn:f_star_omega_11bar} & & (f^\star\omega)\bigg(\frac{\partial}{\partial z_1}_{|(z_1^0,\,z_2^0)},\,\frac{\partial}{\partial\bar{z}_1}_{|(z_1^0,\,z_2^0)}\bigg) \\
  \nonumber & = & \frac{i}{2}\,(\alpha\wedge\bar\alpha + \beta\wedge\bar\beta + \gamma\wedge\bar\gamma)\bigg(-ie^{-z_1^0}z_2^0\,A_{g_0} + e^{-z_1^0}z_2^0\,B_{g_0}-2iC_{g_0},\, ie^{-\overline{z_1^0}}\overline{z_2^0}\,\bar{A}_{g_0} + e^{-\overline{z_1^0}}\overline{z_2^0}\,\bar{B}_{g_0}+2i\bar{C}_{g_0}\bigg) \\
  & = & i\,(|z_2^0|^2\,e^{-2\,\mbox{\small Re}\,(z_1^0)} + 2).\end{eqnarray}

 From (\ref{eqn:f_star_z_2}), we get: \begin{eqnarray}\label{eqn:f_star_omega_22bar} & & (f^\star\omega)\bigg(\frac{\partial}{\partial z_2}_{|(z_1^0,\,z_2^0)},\,\frac{\partial}{\partial\bar{z}_2}_{|(z_1^0,\,z_2^0)}\bigg) \\
   \nonumber & = & \frac{i}{2}\,(\alpha\wedge\bar\alpha + \beta\wedge\bar\beta + \gamma\wedge\bar\gamma)\bigg(-ie^{-z_1^0}\,A_{g_0} + e^{-z_1^0}\,B_{g_0},\, ie^{-\overline{z_1^0}}\,\bar{A}_{g_0} + e^{-\overline{z_1^0}}\,\bar{B}_{g_0}\bigg) = i\,e^{-2\,\mbox{\small Re}\,(z_1^0)}.\end{eqnarray}

 From (\ref{eqn:f_star_z_1}) and (\ref{eqn:f_star_z_2}), we get: \begin{eqnarray}\nonumber\label{eqn:f_star_omega_12bar} & & (f^\star\omega)\bigg(\frac{\partial}{\partial z_1}_{|(z_1^0,\,z_2^0)},\,\frac{\partial}{\partial\bar{z}_2}_{|(z_1^0,\,z_2^0)}\bigg) \\
  \nonumber & = & \frac{i}{2}\,(\alpha\wedge\bar\alpha + \beta\wedge\bar\beta + \gamma\wedge\bar\gamma)\bigg(-ie^{-z_1^0}z_2^0\,A_{g_0} + e^{-z_1^0}z_2^0\,B_{g_0}-2iC_{g_0},\, ie^{-\overline{z_1^0}}\,\bar{A}_{g_0} + e^{-\overline{z_1^0}}\,\bar{B}_{g_0}\bigg) \\
  & = & i\,z_2^0\,e^{-2\,\mbox{\small Re}\,(z_1^0)}.\end{eqnarray}

Finally, from (\ref{eqn:f_star_z_1}) and (\ref{eqn:f_star_z_2}) we also get: \begin{eqnarray}\nonumber\label{eqn:f_star_omega_21bar} & & (f^\star\omega)\bigg(\frac{\partial}{\partial z_2}_{|(z_1^0,\,z_2^0)},\,\frac{\partial}{\partial\bar{z}_1}_{|(z_1^0,\,z_2^0)}\bigg) \\
  \nonumber & = & \frac{i}{2}\,(\alpha\wedge\bar\alpha + \beta\wedge\bar\beta + \gamma\wedge\bar\gamma)\bigg(-ie^{-z_1^0}\,A_{g_0} + e^{-z_1^0}\,B_{g_0},\, ie^{-\overline{z_1^0}}\,\overline{z_2^0}\,\bar{A}_{g_0} + e^{-\overline{z_1^0}}\,\overline{z_2^0}\,\bar{B}_{g_0} +2i\,\overline{C}_{g_0}\bigg) \\
  & = & i\,\overline{z_2^0}\,e^{-2\,\mbox{\small Re}\,(z_1^0)}.\end{eqnarray}

All that is left to do is to put (\ref{eqn:f_star_omega_11bar}), (\ref{eqn:f_star_omega_22bar}), (\ref{eqn:f_star_omega_12bar}) and (\ref{eqn:f_star_omega_21bar}) together and get the contention.  \hfill $\qed$

 \vspace{3ex}

 Based on this, an elementary calculation, spelt out in the proof of the following statement, shows that $f$ is not of subexponential growth.

 \begin{Lem}\label{Lem:SL2C_f-star-omega_computation_1} The map $f$ defined in (\ref{eqn:hol-map_C2-SL2C}) has the following property: \begin{eqnarray*}\log\int\limits_0^b\mbox{Vol}_{\omega,\,f}(B_t)\,dt\geq\sqrt{2}\,b, \hspace{3ex} b\gg 1.\end{eqnarray*}

 \end{Lem}

 \noindent {\it Proof.} Taking squares in the expression for $f^\star\omega$ of Lemma \ref{Lem:SL2C_f-star-omega_computation}, we get in $\C^2\simeq\R^4 $: \begin{eqnarray*}f^\star\omega^2 = 4e^{-2\mbox{Re}\,z_1}\,idz_1\wedge d\bar{z}_1\wedge idz_2\wedge d\bar{z}_2 = 16e^{-2x_1}\,dx_1\wedge dy_1\wedge dx_2\wedge dy_2,\end{eqnarray*} where $z_1=x_1 +iy_1$ and $z_2=x_2 + iy_2$.

 Passing to spherical coordinates $(\rho,\,\theta_1,\,\theta_2,\,\theta_3)$ in $\R^4$, with $\rho\geq 0$, $\theta_1,\theta_2\in[0,\,\pi]$, $\theta_3\in[0,\,2\pi)$, such that $x_1=\rho\cos\theta_1$, $y_1= \rho(\sin\theta_1)(\cos\theta_2)$, $x_2= \rho(\sin\theta_1)(\sin\theta_2)(\cos\theta_3)$ and $y_2= \rho(\sin\theta_1)(\sin\theta_2)(\sin\theta_3)$, we get: \begin{eqnarray*}\mbox{Vol}_{\omega,\,f}(B_t) = \frac{1}{2}\,\int\limits_{B_t}f^\star\omega^2 & = & 8(2\pi^2)\,\int\limits_0^\pi\bigg(\int\limits_0^t e^{-2\rho\cos\theta_1}\,d\rho\bigg)\,d\theta_1  \\
 & \geq & 16\pi^2\,\int\limits_{\frac{3\pi}{4}}^\pi\bigg(\int\limits_0^t e^{-2\rho\cos\theta_1}\,d\rho\bigg)\,d\theta_1  = -8\pi^2\,\int\limits_{\frac{3\pi}{4}}^\pi\frac{e^{-2t\cos\theta_1}}{\cos\theta_1}\,d\theta_1 + a,\end{eqnarray*} where $a\in\R$ is independent of $t$. Since $-1\leq\cos\theta_1\leq -\sqrt{2}/2$ (hence also $1\leq -1/\cos\theta_1\leq\sqrt{2}$) for $\theta_1\in[3\pi/4,\,\pi]$, we get: \begin{eqnarray*}\mbox{Vol}_{\omega,\,f}(B_t) = \frac{1}{2}\,\int\limits_{B_t}f^\star\omega^2\geq 8\pi^2\frac{\pi}{4}\,e^{\sqrt{2}\,t} + a, \hspace{3ex} t>0.\end{eqnarray*} Integrating over $t\in[0,\,b]$, with $b>0$, we get: \begin{eqnarray*}\int\limits_0^b\mbox{Vol}_{\omega,\,f}(B_t)\,dt\geq\frac{2\pi^3}{\sqrt{2}}\,(e^{\sqrt{2}\,b}-1) + a\,b \geq e^{\sqrt{2}\,b}, \hspace{3ex} b\gg 1.\end{eqnarray*}

   This proves the contention.  \hfill $\qed$

   \vspace{2ex}

   We conclude that, for any constant $C>0$, we have: $$\frac{b}{C} - \log\int\limits_0^b\mbox{Vol}_{\omega,\,f}(B_t)\,dt\leq\bigg(\frac{1}{C}-\sqrt{2}\bigg)\,b\longrightarrow -\infty \hspace{3ex}\mbox{as}\hspace{1ex} b\to +\infty$$ if the constant $C$ is chosen such that $C>1/\sqrt{2}$.

   Thus, for any lattice $\Gamma\subset G=SL(2,\,\C)$, the map $f:\C^2\to X=G/\Gamma$ is {\it not of subexponential growth} in the sense of Definition \ref{Def:subexp-growth}.

\vspace{2ex}

(V)\, We now discuss the rather curious example of the $n$-dimensional complex projective space $\Proj^n$ for an arbitrary integer $n\geq 2$. We will see that an obvious map $\C^{n-1}\to\Proj^n$ easily satisfies condition (ii) but may not satisfy condition (i) in the definition \ref{Def:subexp-growth} of the subexponential growth.

\vspace{1ex}


   

 Let $j:\C^{n-1}\to\Proj^n$ be the holomorphic embedding obtained by composing the inclusions $\C^{n-1}\hookrightarrow\C^n$ and $\C^n\hookrightarrow\Proj^n$ given respectively by $(z_1,\dots , \, z_{n-1})\mapsto(z_1,\dots , \, z_{n-1},\,0)$ and $(z_1,\dots , \, z_n)\mapsto[1:z_1:\dots :z_n]$.

As is well known, the restriction to $\C^n$ of the Fubini-Study metric $\omega_{FS}$ of $\Proj^n$ under the above inclusion $\C^n\hookrightarrow\Proj^n$ is \begin{eqnarray*}\omega_{FS} = i\partial\bar\partial\log(1 + |z|^2) = \frac{1}{1 + |z|^2}\,i\partial\bar\partial|z|^2 - \frac{1}{(1 + |z|^2)^2}\,i\partial|z|^2\wedge\bar\partial|z|^2,\end{eqnarray*} where $|z|^2=\sum_{l=1}^n|z_l|^2$. Since $i\partial|z|^2\wedge\bar\partial|z|^2\geq 0$ as a $(1,\,1)$-form and $\omega_0:=(i/2)\,\partial\bar\partial|z|^2$ is the Euclidean metric on $\C^n$, we get: \begin{eqnarray*}\omega_{FS}\leq\frac{2}{1 + |z|^2}\,\omega_0 \leq 2\,\omega_0\hspace{3ex} \mbox{on}\hspace{1ex} \C^n.\end{eqnarray*}

Restricting to $\C^{n-1}$ under the above inclusion $\C^{n-1}\hookrightarrow\C^n$, we get \begin{eqnarray*}j^\star\omega_{FS}\leq 2\beta_0\hspace{3ex} \mbox{on}\hspace{1ex} \C^{n-1},\end{eqnarray*} where $\beta_0=(\omega_0)_{|\C^{n-1}}$ is the Euclidean metric on $\C^{n-1}$. Hence, the $(\omega_{FS},\,j)$-volume of the ball $B_r\subset\C^{n-1}$ of radius $r$ centred at $0$ is estimated as:  \begin{eqnarray*}\mbox{Vol}_{\omega_{FS},\,j}(B_r)=\frac{1}{(n-1)!}\,\int\limits_{B_r}j^\star\omega_{FS}^{n-1}\leq \frac{2^{n-1}}{(n-1)!}\,\int\limits_{B_r}\beta_0^{n-1} = c_n\,r^{2n-2}, \hspace{3ex} r>0, \end{eqnarray*} where $c_n>0$ is a constant depending only on $n$. 

This shows that the embedding $j:\C^{n-1}\to(\Proj^n,\,\omega_{FS})$ is of {\it finite order}, so it satisfies property (ii) in Definition \ref{Def:subexp-growth}. 

\vspace{1ex}

As for property (i) in Definition \ref{Def:subexp-growth}, $j^\star\omega_{FS}$ is a K\"ahler metric on $\C^{n-1}$, hence the second term in formula (\ref{eqn:alternative-formula_integral_growth}) for $\int_{S_t}|d\tau|_{j^\star\omega_{FS}}\,d\sigma_{\omega_{FS},\,j,\,t}$ vanishes. To compute the first term, we deduce from \begin{eqnarray*}j^\star\omega_{FS} = \frac{1}{(1+|z|^2)^2}\,\sum\limits_{j,\,k}\bigg(\delta_{jk}\,(1+|z|^2) - \bar{z}_j\,z_k\bigg)\,idz_j\wedge d\bar{z}_k \hspace{3ex} \mbox{on}\hspace{1ex} \C^{n-1}\end{eqnarray*} that \begin{eqnarray*}\Lambda_{j^\star\omega_{FS}}(i\partial\bar\partial\tau) = (1+|z|^2)^2\,\sum\limits_{j=1}^{n-1}\frac{1}{1+|z|^2 - |z_j|^2}.\end{eqnarray*} Hence, using (\ref{eqn:alternative-formula_integral_growth}) for the equality below, we get: \begin{eqnarray*}\int\limits_{S_t}|d\tau|_{j^\star\omega_{FS}}\,d\sigma_{\omega_{FS},\,j,\,t} & = & 2\,\int\limits_{B_t}(1+|z|^2)^2\,\bigg(\sum\limits_{j=1}^{n-1}\frac{1}{1+|z|^2 - |z_j|^2}\bigg)\,(j^\star\omega_{FS})_{n-1} \\
& \leq & 2(n-1)(1+t^2)^2\,\mbox{Vol}_{\omega_{FS},\,j}(B_t),  \hspace{3ex} t>0.\end{eqnarray*} The last inequality falls far short of the required (\ref{eqn:comparative-growth}), so we cannot say anything at this stage about whether $\Proj^n$ is divisorially hyperbolic or not.

\vspace{2ex}

(VI)\, (a)\, A prototypical example of a compact complex manifold that is {\it not divisorially hyperbolic} is any {\it complex torus} $X=\C^n/\Gamma$, where $\Gamma\subset(\C^n,\,+)$ is any lattice. Any Hermitian metric with constant coefficients on $\C^n$ (for example, the Euclidean metric $\beta = (1/2)\,\sum_j idz_j\wedge d\bar{z}_j$) defines a K\"ahler metric $\omega$ on $X$: $\pi^\star\omega = \beta$, where $\pi:\C^n\to X$ is the projection. If $j:\C^{n-1}\longrightarrow\C^n$ is the obvious inclusion $(z_1,\dots , , z_{n-1})\mapsto(z_1,\dots , , z_n)$, the non-degenerate holomorphic map $f=\pi\circ j: \C^{n-1}\to X$ has subexponential growth thanks to Lemma \ref{Lem:subexp-growth_standard-metric} because $f^\star\omega = j^\star\beta = \beta_0$, where $\beta_0$ is the Euclidean metric of $\C^{n-1}$. This shows that the complex torus $X=\C^n/\Gamma$ is not divisorially hyperbolic.

\vspace{2ex}

(b)\, Similarly, the {\it Iwasawa manifold} $X=G/\Gamma$ is {\it not divisorially hyperbolic}, where $G=(\C^3,\,\star)$ is the nilpotent complex Lie group (called the Heisenberg group) whose group operation is defined as $$(\zeta_1,\,\zeta_2,\,\zeta_3)\star(z_1,\,z_2,\,z_3)=(\zeta_1+z_1,\,\zeta_2+z_2,\,\zeta_3+z_3+\zeta_1\,z_2),$$ while the lattice $\Gamma\subset G$ consists of the elements $ (z_1,\,z_2,\,z_3)\in G$ with $z_1,\,z_2,\,z_3\in\Z[i]$. (See e.g. [Nak75].)

Indeed, the holomorphic $(1,\,0)$-forms $dz_1,\,dz_2,\, dz_3-z_1\,dz_2$ on $\C^3$ induce holomorphic $(1,\,0)$-forms $\alpha,\,\beta,\,\gamma$ on $X$. The Hermitian metric $$\omega_0 = i\alpha\wedge\bar\alpha + i\beta\wedge\bar\beta + i\gamma\wedge\bar\gamma$$ on $X$ lifts to the Hermitian metric $$\omega=\pi^\star\omega_0 = idz_1\wedge d\bar{z}_1 + (1+|z_1|^2)\,idz_2\wedge d\bar{z}_2 + idz_3\wedge d\bar{z}_3 - \bar{z}_1\,idz_3\wedge d\bar{z}_2 - z_1\,idz_2\wedge d\bar{z}_3$$ on $G=\C^3$, where $\pi:G\to X$ is the projection.

Considering the non-degenerate holomorphic map $f=\pi\circ j:\C^2\longrightarrow X$, where $j:\C^2\longrightarrow\C^3$ is the obvious inclusion $(z_1,\, z_2)\mapsto(z_1,\, z_2,\,0)$, we get $$f^\star\omega_0 = j^\star\omega = \omega_{|\C^2} = idz_1\wedge d\bar{z}_1 + (1+|z_1|^2)\,idz_2\wedge d\bar{z}_2$$ on $\C^2$. Hence, $$f^\star\omega_0^2 = 2(1+|z_1|^2)\,dV_0$$ on $\C^2$, where we put $dV_0:=idz_1\wedge d\bar{z}_1\wedge idz_2\wedge d\bar{z}_2$.. Thus, for the ball $B_r\subset\C^2$ of radius $r$ centred at $0$, we get \begin{eqnarray}\label{eqn:Iwasawa-examp_vol-formula}\mbox{Vol}_{\omega_0,\,f}(B_r) = \frac{1}{2}\,\int\limits_{B_r}f^\star\omega_0^2 = \int\limits_{B_r}(1+|z_1|^2)\,dV_0 \leq c_2\,r^4(1+r^2),  \hspace{3ex} r>0,\end{eqnarray} where $c_2>0$ is a constant independent of $r$. This shows that $f$ is of {\it finite order}, hence $f$ satisfies property (ii) in Definition \ref{Def:subexp-growth}.

To show that $f$ has subexponential growth, it remains to check that it also satisfies property (i) in Definition \ref{Def:subexp-growth}. We will first compute the integral on the left of (\ref{eqn:comparative-growth}) in this case. Recall that $n=3$. Then, note that $$d(f^\star\omega_0) = \partial|z_1|^2\wedge idz_2\wedge d\bar{z}_2 + \bar\partial|z_1|^2\wedge idz_2\wedge d\bar{z}_2 = (\partial\tau + \bar\partial\tau)\wedge idz_2\wedge d\bar{z}_2.$$ Thus, we get the following equalities on $\C^2$: \begin{eqnarray*}i(\bar\partial\tau - \partial\tau)\wedge d(f^\star\omega_0) & = & i(\bar\partial\tau - \partial\tau)\wedge(\partial\tau + \bar\partial\tau)\wedge idz_2\wedge d\bar{z}_2 = -2i\partial\tau\wedge\bar\partial\tau\wedge idz_2\wedge d\bar{z}_2 \\
  & = & -2|z_1|^2\,idz_1\wedge d\bar{z}_1\wedge idz_2\wedge d\bar{z}_2= -2|z_1|^2\,dV_0,\end{eqnarray*} where we put $dV_0:=idz_1\wedge d\bar{z}_1\wedge idz_2\wedge d\bar{z}_2$.

On the other hand, since $i\partial\bar\partial\tau = idz_1\wedge d\bar{z}_1 + idz_2\wedge d\bar{z}_2$, we have \begin{eqnarray*}\Lambda_{f^\star\omega_0}(i\partial\bar\partial\tau) = 1 + \frac{1}{1+|z_1|^2}.\end{eqnarray*}

Therefore, the integral on the left of (\ref{eqn:comparative-growth}) reads in this case: \begin{eqnarray*}\int\limits_{S_t}|d\tau|_{f^\star\omega_0}\,d\sigma_{\omega_0,\,f,\,t} & = & 2\int\limits_{B_t}\bigg(1 + \frac{1}{1+|z_1|^2}\bigg)\,(1 + |z_1|^2)\,dV_0 + 2\int\limits_{B_t}|z_1|^2\,dV_0 \\
  & = & 4\,\int\limits_{B_t}(1+|z_1|^2)\,dV_0 = 4\,\mbox{Vol}_{\omega_0,\,f}(B_t), \hspace{3ex} t>0,\end{eqnarray*} where the last equality was seen in (\ref{eqn:Iwasawa-examp_vol-formula}).

This proves that $f$ satisfies property (i) in Definition \ref{Def:subexp-growth}. We conclude that the map $f$ has {\it subexponential growth}, proving that the $3$-dimensional Iwasawa manifold $X$ is not divisorially hyperbolic.

\vspace{2ex}

\vspace{2ex}

(c)\, Finally, we point out that no {\it Nakamura manifold} $X=G/\Gamma$ is divisorially hyperbolic, where $G=(\C^3,\,\star)$ is the solvable, non-nilpotent complex Lie group whose group operation is defined as $$(\zeta_1,\,\zeta_2,\,\zeta_3)\star(z_1,\,z_2,\,z_3)=(\zeta_1+z_1,\,\zeta_2+e^{-\zeta_1}z_2,\,\zeta_3+e^{\zeta_1}z_3),$$ while $\Gamma\subset G$ is a lattice. (See e.g. [Nak75].)

We equip $X$ with the metric $$\omega_0 = i\eta_1\wedge\bar\eta_1 + i\eta_2\wedge\bar\eta_2 + i\eta_3\wedge\bar\eta_3,$$ where $\eta_1,\eta_2,\eta_3$ are the holomorphic $(1,\,0)$-forms on $X$ induced respectively by the left-invariant holomorphic $(1,\,0)$-forms $dz_1,\,e^{-z_1}dz_2,\,e^{z_1}dz_3$ on $G$. If $\pi:G\longrightarrow X$ is the projection, we see that $\omega_0$ lifts to the Hermitian metric $$\omega=\pi^\star\omega_0 = idz_1\wedge d\bar{z}_1 + e^{-2\mbox{Re}(z_1)}\,idz_2\wedge d\bar{z}_2 + e^{2\mbox{Re}(z_1)}\,idz_3\wedge d\bar{z}_3$$ on $G=\C^3$.

Let $j:\C^2\longrightarrow G\simeq\C^3$ be the obvious inclusion $(z_2,\,z_3)\mapsto(0,\,z_2,\,z_3)$. Then, for the non-degenerate holomorphic map $f=\pi\circ j:\C^2\longrightarrow X$, we get $$f^\star\omega_0 = j^\star(\pi^\star\omega_0) = idz_2\wedge d\bar{z}_2 + idz_3\wedge d\bar{z}_3$$ on $\C^2$. Thus, $f^\star\omega_0$ is the Euclidean metric of $\C^2$, so $f$ has {\it subexponential growth} by Lemma \ref{Lem:subexp-growth_standard-metric}, proving that the Nakamura manifold $X$ is {\it not divisorially hyperbolic}.

\section{Divisorially K\"ahler and divisorially nef classes}\label{section:div-K_div-nef_classes}

The starting point is the following simple, but key observation.

\begin{Lem}\label{Lem:P_map} Let $X$ be a compact complex manifold with $\mbox{dim}_\C X=n$. The map: \begin{equation}\label{eqn:P_n-1_map_def}P=P^{n-1}_{n-1,\,n-1}:H^2_{DR}(X,\,\R)\longrightarrow  H^{n-1,\,n-1}_A(X,\,\R), \hspace{3ex} \{\alpha\}_{DR}\longmapsto\{(\alpha^{n-1})^{n-1,\,n-1}\}_A,\end{equation} is {\bf well defined} in the sense that it is independent of the choice of a $C^\infty$ representative $\alpha$ of its De Rham cohomology class, where $(\alpha^{n-1})^{n-1,\,n-1}$ is the component of bidegree $(n-1,\,n-1)$ of the $(2n-2)$-form $\alpha^{n-1}$.

\end{Lem}  

This follows from the following

\begin{Lem}\label{Lem:correctness_div-nef_classes} Let $X$ be a compact complex manifold with $\mbox{dim}_\C X=n$.

  \vspace{1ex}

  (i)\, For any $k\in\{0,\dots , 2n\}$, any form $\alpha\in C^\infty_k(X,\,\C)$ such that $d\alpha=0$ and any bidegree $(p,\,q)$ with $p+q=k$, we have $$\partial\bar\partial\alpha^{p,\,q}=0,$$ where $\alpha^{p,\,q}$ is the $(p,\,q)$-type component of $\alpha$.

  In particular, for every $2$-form $\alpha$ such that $d\alpha=0$, we have $\partial\bar\partial(\alpha^{n-1})^{n-1,\,n-1}=0$, so $(\alpha^{n-1})^{n-1,\,n-1}$ defines an Aeppli cohomology class $\{(\alpha^{n-1})^{n-1,\,n-1}\}_A\in H^{n-1,\,n-1}_A(X,\,\C)$.

  \vspace{1ex}

  (ii)\, For any $2$-forms $\alpha_1$ and $\alpha_2$ such that $d\alpha_1= d\alpha_2=0$ and $\alpha_1 = \alpha_2 + d\beta$ for some $1$-form $\beta$, we have $$\{(\alpha_1^{n-1})^{n-1,\,n-1}\}_A  = \{(\alpha_2^{n-1})^{n-1,\,n-1}\}_A.$$

\end{Lem}

\noindent {\it Proof.} (i)\, Writing the decomposition $\alpha = \sum\limits_{r+s=k}\alpha^{r,\,s}$ of $\alpha$ into pure-type forms, we see that the hypothesis $d\alpha=0$ is equivalent to $\partial\alpha^{r,\,s} + \bar\partial\alpha^{r+1,\,s-1} = 0$ for all $(r,\,s)$. Applying $\bar\partial$, this implies that $\partial\bar\partial\alpha^{r,\,s}=0$ for all $(r,\,s)$.

(ii)\, Taking the $(n-1)$-st power in $\alpha_1 = \alpha_2 + d\beta$ and using the fact that $d\alpha_2=0$, we get: $$\alpha_1^{n-1} = \alpha_2^{n-1} + \sum\limits_{k=1}^{n-1}{n-1 \choose k}\,d(\alpha_2^{n-1-k}\wedge\beta\wedge(d\beta)^{k-1}).$$ Hence, $\alpha_1^{n-1} - \alpha_2^{n-1}\in\mbox{Im}\,d$, which implies that $(\alpha_1^{n-1})^{n-1,\,n-1}-(\alpha_2^{n-1})^{n-1,\,n-1}\in\mbox{Im}\,\partial + \mbox{Im}\,\bar\partial$. \hfill $\qed$

\subsection{Case of projective manifolds}\label{subsection:div-nef_classes_proj-manif}

Let $X$ be a {\bf projective} manifold with $\mbox{dim}_\C X=n$. As is well known (see e. g. [Ha70]), a holomorphic line bundle $L\longrightarrow X$ is said to be {\bf nef} if $$L.C:=\int\limits_C c_1(L)\geq 0$$ for every {\it curve} $C\subset X$, where $c_1(L)=\{\frac{i}{2\pi}\Theta_h(L)\}_{DR}\in H^{1,\,1}(X,\,\R)\cap H^2(X,\,\Z)$ is the first Chern class of $L$, namely the De Rham cohomology class of the curvature form of $L$ with respect to any Hermitian metric $h$ on $L$.

\vspace{2ex}

We now generalise this notion in the context of divisors (rather than curves) and of possibly non-integral and non-type $(1,\,1)$ real De Rham cohomology classes using the Serre-type duality (\ref{eqn:BC-A_duality}).  

\begin{Def}\label{Def:div-nef_classes_proj-manif} Let $X$ be a {\bf projective} manifold with $\mbox{dim}_\C X=n$. A cohomology class $\{\alpha\}_{DR}\in H^2_{DR}(X,\,\R)$ is said to be {\bf projectively divisorially nef} if $$P(\{\alpha\}_{DR}).\{[D]\}_{BC}:=\int\limits_D(\alpha^{n-1})^{n-1,\,n-1}\geq 0$$ for all effective {\bf divisors} $D\geq 0$ on $X$ and some (hence any) representative $\alpha\in C^\infty_2(X,\,\R)$ of $\{\alpha\}_{DR}$.

\end{Def}

\vspace{3ex}

As is well known, the current of integration $[D]$ on an effective divisor $D$ is a closed positive $(1,\,1)$-current. (By a $(1,\,1)$-current we mean a current of {\it bidegree} $(1,\,1)$.) However, not every such current $T$ is the current of integration on an effective divisor $D$. Nevertheless, we notice that Definition \ref{Def:div-nef_classes_proj-manif} does not change if divisors are replaced by currents whose cohomology classes lie in the real vector space $NS_\R(X)$ spanned by {\it integral} $(1,\,1)$-cohomology classes, known as the Neron-Severi subspace of $H^2_{DR}(X,\,\R)$. 

\begin{Prop}\label{Prop:div-nef_classes_currents} Let $X$ be a {\bf projective} manifold with $\mbox{dim}_\C X=n$. Fix any cohomology class $\{\alpha\}_{DR}\in H^2_{DR}(X,\,\R)$. Then, $\{\alpha\}_{DR}$ is {\bf projectively divisorially nef} if and only if $$\{(\alpha^{n-1})^{n-1,\,n-1}\}_A.\{T\}_{BC}:=\int\limits_X(\alpha^{n-1})^{n-1,\,n-1}\wedge T\geq 0$$ for all closed positive $(1,\,1)$-currents $T\geq 0$ on $X$ such that $\{T\}_{DR}\in NS_\R(X)$ and some (hence any) representative $\alpha\in C^\infty_2(X,\,\R)$ of $\{\alpha\}_{DR}$.

\end{Prop}

\noindent {\it Proof.} It is proved in (b) of Proposition 6.6. in [Dem00], as a consequence of Nadel's Vanishing Theorem, that the integral part of the  pseudo-effective cone of $X$ (i.e. the set of cohomology classes of closed positive $(1,\,1)$-currents that are linear combinations with real coefficients of integral  classes) is the closure of the {\it effective cone} of $X$ (i.e. the set of cohomology classes of effective divisors). This means that, for every closed positive $(1,\,1)$-current $T\geq 0$ on $X$ such that $\{T\}_{DR}\in NS_\R(X)$, the class $\{T\}_{BC}$ is a limit of classes $\{[D_j]\}_{BC}$ with $(D_j)_{j\in\N}$ effective divisors on $X$. This suffices to conclude. \hfill $\qed$ 

\vspace{2ex}

We now observe a useful property of the set of projectively divisorially nef classes.

\begin{Prop}\label{Prop:div-nef_cone_closed} Let $X$ be a {\bf projective} manifold. The set $${\cal PDN}_X:=\bigg\{\{\alpha\}_{DR}\in H^2_{DR}(X,\,\R)\,\mid\, \{\alpha\}_{DR} \hspace{1ex}\mbox{is projectively divisorially nef}\hspace{1ex}\bigg\}$$ is a {\bf closed} cone in $H^2_{DR}(X,\,\R)$.

\end{Prop}

\noindent {\it Proof.} To show that ${\cal PDN}_X$ is a cone, we have to show that it is stable under multiplications by non-negative reals, which is obvious.

To show closedness, let $\{\alpha\}_{DR}\in H^2_{DR}(X,\,\R)$ be a limit of classes $\{\alpha_j\}_{DR}\in{\cal PDN}_X$. Then, the class $\{\alpha_j\}_{DR} - \{\alpha\}_{DR}$ converges to $0$ in $H^2_{DR}(X,\,\R)$ as $j\to +\infty$. By the definition of the quotient topology of $H^2_{DR}(X,\,\R)$, there exists a sequence of $C^\infty$ $d$-closed $2$-forms $\beta_j\in\{\alpha_j\}_{DR} - \{\alpha\}_{DR}$ such that $\beta_j\longrightarrow 0$ in the $C^\infty$ topology as $j\to +\infty$. Now, pick an arbitrary $C^\infty$ representative $\alpha$ of the class $\{\alpha\}_{DR}$. We infer that $\alpha_j:=\alpha + \beta_j$ represents the class $\{\alpha_j\}_{DR}$ for every $j$ and $\lim_{j\to +\infty}\alpha_j = \alpha$ in the $C^\infty$ topology.

 Thus, for every effective divisor $D$ on $X$, we have: $$\int\limits_D(\alpha_j^{n-1})^{n-1,\,n-1}\geq 0 \hspace{2ex} \mbox{for all} \hspace{2ex} j\in\N \hspace{2ex} \mbox{and} \hspace{2ex} \lim_{j\to +\infty}\int\limits_D(\alpha_j^{n-1})^{n-1,\,n-1} = \int\limits_D(\alpha^{n-1})^{n-1,\,n-1},$$ where the first inequality follows from the assumption $\{\alpha_j\}_{DR}\in{\cal PDN}_X$ for all $j$. Therefore, $\int_D(\alpha^{n-1})^{n-1,\,n-1}\geq 0$ for every effective divisor $D$, proving that the class $\{\alpha\}_{DR}$ is projectively divisorially nef.  \hfill $\qed$

\subsection{Case of arbitrary compact complex manifolds}\label{subsection:div-nef_classes_arbitrary-manif}

Let $X$ be a compact complex $n$-dimensional manifold. Recall the obvious inclusion ${\cal SG}_X\subset{\cal G}_X$ of the {\it strongly Gauduchon (sG) cone} of $X$ in the {\it Gauduchon cone}. (See the introduction for a reminder of the definitions.) The equality ${\cal SG}_X = {\cal G}_X$ is equivalent to every Gauduchon metric on $X$ being sG (see [PU18, Lemma 1.3]). Manifolds $X$ with this property are called {\bf sGG manifolds}; they were studied in [PU18].

The definition of divisorially nef classes given in (ii) of the next Definition \ref{Def:div-K-nef_classes} on an arbitrary compact complex manifold will be shown in (b) of Proposition \ref{Prop:div-nef-cone_prop} to imply the projectively divisorially nef property defined on projective manifolds in Definition \ref{Def:div-nef_classes_proj-manif}.

\begin{Def}\label{Def:div-K-nef_classes} Let $X$ be a compact complex manifold with $\mbox{dim}_\C X=n$ and let \begin{equation*}P:H^2_{DR}(X,\,\R)\longrightarrow  H^{n-1,\,n-1}_A(X,\,\R), \hspace{3ex} \{\alpha\}_{DR}\longmapsto\{(\alpha^{n-1})^{n-1,\,n-1}\}_A,\end{equation*} be the map of Lemma \ref{Lem:P_map}.

  \vspace{1ex}

  (i)\, A cohomology class $\{\alpha\}_{DR}\in H^2_{DR}(X,\,\R)$ is said to be {\bf divisorially K\"ahler} if $P(\{\alpha\}_{DR})\in{\cal G}_X$. The set $${\cal DK}_X:=\bigg\{\{\alpha\}_{DR}\in H^2_{DR}(X,\,\R)\,\mid\, \{\alpha\}_{DR} \hspace{1ex}\mbox{is divisorially K\"ahler}\hspace{1ex}\bigg\}$$ is called the {\bf divisorially K\"ahler cone} of $X$. 

  \vspace{1ex}

  (ii)\, A cohomology class $\{\alpha\}_{DR}\in H^2_{DR}(X,\,\R)$ is said to be {\bf divisorially nef} if $P(\{\alpha\}_{DR})\in\overline{\cal G}_X$, where $\overline{\cal G}_X$ is the closure of the Gauduchon cone in $H^{n-1,\,n-1}_A(X,\,\R)$.

  The set $${\cal DN}_X:=\bigg\{\{\alpha\}_{DR}\in H^2_{DR}(X,\,\R)\,\mid\, \{\alpha\}_{DR} \hspace{1ex}\mbox{is divisorially nef}\hspace{1ex}\bigg\}$$ is called the {\bf divisorially nef cone} of $X$.

\end{Def}

Note that ${\cal DK}_X$ and ${\cal DN}_X$ are {\it cones} in $H^2_{DR}(X,\,\R)$ in the sense that they are stable under multiplications by positive reals. However, they are not convex and are not stable under additions since the map $P$ is not linear.

\begin{Prop}\label{Prop:div-K-cone_prop} Let $X$ be a compact complex manifold with $\mbox{dim}_\C X = n$. The {\bf divisorially K\"ahler cone} of $X$ can be described as \begin{equation}\label{eqn:div-K-cone_description}{\cal DK}_X = P^{-1}({\cal SG}_X) = P^{-1}({\cal G}_X)\subset H^2_{DR}(X,\,\R).\end{equation}

  In particular, ${\cal DK}_X$ is {\bf open} in $H^2_{DR}(X,\,\R)$ and the following implication holds: \begin{equation}\label{eqn:div-K-cone-sG_implication}{\cal DK}_X\neq\emptyset\implies X \hspace{1ex}\mbox{is an sG manifold}.\end{equation}

\end{Prop}

\noindent {\it Proof.} The identity ${\cal DK}_X = P^{-1}({\cal G}_X)$ holds by the definition of ${\cal DK}_X$. Meanwhile, $P^{-1}({\cal SG}_X)\subset P^{-1}({\cal G}_X)$ since ${\cal SG}_X\subset{\cal G}_X$. So, it suffices to prove the inclusion ${\cal DK}_X\subset P^{-1}({\cal SG}_X)$.

Let $\{\alpha\}_{DR}\in{\cal DK}_X$. Pick an arbitrary smooth representative $\alpha\in\{\alpha\}_{DR}$. Since $P(\{\alpha\}_{DR})\in{\cal G}_X$, there exists a Gauduchon metric $\omega$ on $X$ such that $$(\alpha^{n-1})^{n-1,\,n-1} = \omega^{n-1} + \partial u^{n-2,\,n-1} + \bar\partial u^{n-1,\,n-2}$$ for some smooth forms $u^{n-2,\,n-1}$ and $u^{n-1,\,n-2}$ of the displayed bidegrees. These forms can be chosen to be conjugate to each other since $\alpha$ and $\omega$ are real. We get: $$(\alpha^{n-1})^{n-1,\,n-1} = \omega^{n-1} + (d(u^{n-2,\,n-1} + u^{n-1,\,n-2}))^{n-1,\,n-1},$$ so $\omega^{n-1}$ is the $(n-1,\,n-1)$-component of the smooth real $d$-closed $(2n-2)$-form $\alpha^{n-1} - d(u^{n-2,\,n-1} + u^{n-1,\,n-2})$. This proves that $\omega$ is strongly Gauduchon (see [Pop13, Proposition 4.2.]). Since $\{\omega^{n-1}\}_A = \{(\alpha^{n-1})^{n-1,\,n-1}\}_A = P(\{\alpha\}_{DR})$, we infer that $P(\{\alpha\}_{DR})\in{\cal SG}_X$. This proves the inclusion ${\cal DK}_X\subset P^{-1}({\cal SG}_X)$.

The openness of ${\cal DK}_X$ in $H^2_{DR}(X,\,\R)$ follows from the openness of 
${\cal G}_X$ in $H^{n-1,\,n-1}_A(X,\,\R)$ and from the continuity of the map $P$.

Finally, to prove implication (\ref{eqn:div-K-cone-sG_implication}), suppose there exists $\{\alpha\}_{DR}\in{\cal DK}_X$. Then $P(\{\alpha\}_{DR})\in{\cal SG}_X$, so ${\cal SG}_X\neq\emptyset$. The last piece of information is equivalent to $X$ being an sG manifold.  \hfill $\qed$

\vspace{2ex}

Part (c) of the following result gives, on any compact complex manifold, an alternative definition of a {\it divisorially nef} class that is analogous to the classical analytic definition of a nef class given in [Dem92].

\begin{Prop}\label{Prop:div-nef-cone_prop} Let $X$ be a compact complex manifold with $\mbox{dim}_\C X = n$.

  \vspace{1ex}

  (a)\, The divisorially nef cone ${\cal DN}_X$ is {\bf closed} in $H^2_{DR}(X,\,\R)$. In particular, \begin{equation}\label{eqn:div-K-nef_cones_inclusion}\overline{\cal DK}_X\subset{\cal DN}_X,\end{equation} where $\overline{\cal DK}_X$ is the closure of the divisorially K\"ahler cone in $H^2_{DR}(X,\,\R)$.

  \vspace{1ex}

  (b)\, For every class $\{\alpha\}_{DR}\in H^2_{DR}(X,\,\R)$, the following equivalence holds: \begin{equation}\label{eqn:div-nef-cone_divisor_equiv}\{\alpha\}_{DR}\in{\cal DN}_X \iff P(\{\alpha\}_{DR}).\{T\}_{BC}\geq 0 \hspace{3ex}\mbox{for every} \hspace{1ex} \{T\}_{BC}\in{\cal E}_X,\end{equation} where ${\cal E}_X\subset H^{1,\,1}_{BC}(X,\,\R)$ is the pseudo-effective cone of $X$ consisting of the Bott-Chern cohomology classes of all closed positive $(1,\,1)$-currents $T\geq 0$ on $X$.

  In particular, if $X$ is projective, a class $\{\alpha\}_{DR}\in H^2_{DR}(X,\,\R)$ is projectively divisorially nef in the sense of Definition \ref{Def:div-nef_classes_proj-manif} whenever $\{\alpha\}_{DR}\in H^2_{DR}(X,\,\R)$ is divisorially nef in the sense of  Definition \ref{Def:div-K-nef_classes}.

  \vspace{1ex}

  (c)\, A class $\{\alpha\}_{DR}\in H^2_{DR}(X,\,\R)$ is {\bf divisorially nef} if and only if for every constant $\varepsilon>0$, there exists a representative $\Omega_\varepsilon\in C^\infty_{n-1,\,n-1}(X,\,\R)$ of the class $P(\{\alpha\}_{DR})$ such that $$\Omega_\varepsilon\geq -\varepsilon\,\omega^{n-1},$$ where $\omega>0$ is an arbitrary Hermitian metric on $X$ fixed beforehand.

\end{Prop}

\noindent {\it Proof.} Since $X$ is compact, any two Hermitian metrics on $X$ are comparable. Meanwhile, a Gauduchon metric always exists on $X$ by [Gau77a], so we may assume that the background metric $\omega$ on $X$ is actually Gauduchon.

\vspace{1ex}

(a)\, The first statement is an immediate consequence of the continuity of $P$ and of the identity ${\cal DN}_X = P^{-1}(\overline{{\cal G}}_X)$ defining the divisorially nef cone. Inclusion (\ref{eqn:div-K-nef_cones_inclusion}) follows from the closedness of ${\cal DN}_X$ and from the obvious inclusion ${\cal DK}_X\subset{\cal DN}_X$. 

\vspace{1ex}

(b)\, The first statement follows from the duality between the pseudo-effective cone ${\cal E}_X\subset H^{1,\,1}_{BC}(X,\,\R)$ and the closure of the Gauduchon cone $\overline{\cal G}_X\subset H^{n-1,\,n-1}_A(X,\,\R)$ under duality (\ref{eqn:BC-A_duality}) between $H^{1,\,1}_{BC}(X,\,\C)$ and $H^{n-1,\,n-1}_A(X,\,\C)$. This cone duality, observed in [Pop15b] as a reformulation of Lamari's duality Lemma 3.3 in [Lam99], implies that, given any class $\mathfrak{c}^{n-1,\,n-1}_A\in H^{n-1,\,n-1}_A(X,\,\R)$, the following equivalence holds: $$\mathfrak{c}^{n-1,\,n-1}_A\in\overline{\cal G}_X\iff\mathfrak{c}^{1,\,1}_{BC}.\mathfrak{c}^{n-1,\,n-1}_A\geq 0 \hspace{2ex} \mbox{for every class} \hspace{1ex} \mathfrak{c}^{1,\,1}_{BC}\in{\cal E}_X.$$

In our case, it suffices to apply this duality to $\mathfrak{c}^{n-1,\,n-1}_A:=P(\{\alpha\}_{DR})$ to get equivalence (\ref{eqn:div-nef-cone_divisor_equiv}).

When $X$ is projective, the second statement follows from (\ref{eqn:div-nef-cone_divisor_equiv}) and from Proposition \ref{Prop:div-nef_classes_currents}. 

\vspace{1ex}

(c)\, ``$\Longleftarrow$'' Fix a {\it Gauduchon} metric $\omega$ on $X$ and a class $\{\alpha\}_{DR}\in H^2_{DR}(X,\,\R)$. Suppose that, for every $\varepsilon>0$, the class $P(\{\alpha\}_{DR})$ can be represented by a form $\Omega_\varepsilon\in C^\infty_{n-1,\,n-1}(X,\,\R)$ such that $\Omega_\varepsilon\geq -\varepsilon\,\omega^{n-1}$. Then, $\Omega_\varepsilon + 2\varepsilon\,\omega^{n-1} \geq \varepsilon\,\omega^{n-1} > 0$ and $\Omega_\varepsilon + 2\varepsilon\,\omega^{n-1}$ is $\partial\bar\partial$-closed, so it is the $(n-1)$-st power of a Gauduchon metric. (See [Mic83] for the existence of a unique $(n-1)$-st root for any positive definite $(n-1,\,n-1)$-form on an $n$-dimensional complex manifold.) Hence, $$\mathfrak{c}_\varepsilon:=\{\Omega_\varepsilon\}_A + 2\varepsilon\,\{\omega^{n-1}\}_A = P(\{\alpha\}_{DR}) + 2\varepsilon\,\{\omega^{n-1}\}_A\in{\cal G}_X, \hspace{6ex} \varepsilon>0,$$ so $\mathfrak{c}_\varepsilon\longrightarrow P(\{\alpha\}_{DR})$ as $\varepsilon\to 0$. Therefore, $P(\{\alpha\}_{DR})\in\overline{\cal G}_X$, which amounts to $\{\alpha\}_{DR}$ being divisorially nef.

\vspace{1ex}

``$\Longrightarrow$'' Suppose $\{\alpha\}_{DR}\in H^2_{DR}(X,\,\R)$ is divisorially nef. Then $P(\{\alpha\}_{DR})\in\overline{\cal G}_X$, so there exists a sequence of classes $\mathfrak{c}_k\in{\cal G}_X$ such that $\mathfrak{c}_k\longrightarrow P(\{\alpha\}_{DR})$ as $k\to +\infty$. Thus, $P(\{\alpha\}_{DR}) - \mathfrak{c}_k\longrightarrow 0$ in $H^{n-1,\,n-1}_A(X,\,\R)$, so, from the definition of the quotient topology, we infer the existence of a sequence of real $C^\infty$ representatives $\Gamma_k\in P(\{\alpha\}_{DR}) - \mathfrak{c}_k$ such that $\Gamma_k\longrightarrow 0$ in the $C^\infty$ topology (hence also in the $C^0$ topology) as $k\to +\infty$. This implies that, for every $\varepsilon>0$, there exists $k_\varepsilon\in\N$ such that $$\Gamma_k\geq -\varepsilon\,\omega^{n-1}, \hspace{6ex} k\geq k_\varepsilon,$$ where $\omega$ is an arbitrarily fixed Gauduchon metric on $X$.

On the other hand, for every $k\in\N$, pick a Gauduchon metric $\omega_k$ on $X$ such that $\omega_k^{n-1}\in\mathfrak{c}_k$. (This is possible since $\mathfrak{c}_k\in{\cal G}_X$.) We infer that $$P(\{\alpha\}_{DR})\ni\Omega_\varepsilon : = \Gamma_{k_\varepsilon} + \omega_{k_\varepsilon}^{n-1}\geq -\varepsilon\,\omega^{n-1},  \hspace{6ex} \varepsilon>0.$$ This proves the contention.  \hfill $\qed$

\begin{Question}\label{Question:div-K-nef_cones_equality} If ${\cal DK}_X\neq\emptyset$, is (\ref{eqn:div-K-nef_cones_inclusion}) an equality?

\end{Question}

\subsection{Examples}\label{subsection:div-nef_examp} (1)\, Let $X$ be a compact complex manifold with $\mbox{dim}_\C X=n$. Suppose there exists a {\it Hermitian-symplectic (H-S) structure} on $X$. According to [ST10, Definition 1.5], this is a real $C^\infty$ $d$-closed $2$-form $\widetilde\omega$ on $X$ whose component $\omega$ of type $(1,\,1)$ is positive definite. It is easy to see that the $(n-1,\,n-1)$-component $(\widetilde\omega^{n-1})^{n-1,\,n-1}$ of the $(n-1)$-st power of $\widetilde\omega$ is {\it positive definite} on $X$. (See [YZZ19, $\S.2$, Lemma $1$] or [DP20, Proposition 2.1].) Since it is also $\partial\bar\partial$-closed, it defines an element $P(\{\widetilde\omega\}_{DR})=\{(\widetilde\omega^{n-1})^{n-1,\,n-1}\}_A$ in the Gauduchon cone ${\cal G}_X$. Thus, we get

\begin{Prop}\label{Prop:H-S_div-K} If $\widetilde\omega$ is a {\bf Hermitian-symplectic} structure on a compact complex manifold $X$, the cohomology class $\{\widetilde\omega\}_{DR}\in H^2_{DR}(X,\,\R)$ is {\bf divisorially K\"ahler}.

\end{Prop}

(2)\, From the equivalence of (i) and (iii) in Proposition \ref{Prop:deg-bal_characterisation}, we deduce the following

\begin{Prop}\label{Prop:deg-bal_div-K} If $X$ is a {\bf degenerate balanced} compact complex manifold, ${\cal DK}_X = H^2_{DR}(X,\,\R)$.

\end{Prop}

\vspace{2ex}

(3)\, Let $L\longrightarrow X$ be a $C^\infty$ (not necessarily holomorphic) complex line bundle over an $n$-dimensional compact complex manifold $X$. For any $C^\infty$ Hermitian metric $h$ on $L$, the curvature form $\alpha=\frac{i}{2\pi}\Theta_h(L)$ is a $C^\infty$ real $d$-closed $2$-form on $X$ which represents the first Chern class $c_1(L)\in H^2(X,\,\Z)$ of $X$. Moreover, $c_1(L)$ is of type $(1,\,1)$ if and only if $L$ is {\it holomorphic}. We say that $L$ is divisorially nef if $c_1(L)$ is.

Similarly, when $X$ is {\it projective}, we say that $L$ is projectively divisorially nef if $c_1(L)$ is. In this case, we have:

\vspace{1ex}

$\displaystyle L \hspace{2ex} \mbox{is {\bf projectively divisorially nef}} \hspace{2ex} \iff $

\hspace{25ex} $\displaystyle P(c_1(L)).\{[D]\}_{BC} :=\int\limits_D\bigg(\bigg(\frac{i}{2\pi}\Theta_h(L)\bigg)^{n-1}\bigg)^{n-1,\,n-1}\geq 0$

\vspace{1ex}

\noindent for all effective divisors $D$ on $X$, where $P:H^2_{DR}(X,\,\R)\longrightarrow  H^{n-1,\,n-1}_A(X,\,\R)$ is the map (\ref{eqn:P_n-1_map_def}).

  \vspace{2ex}

  (4)\, If $L\longrightarrow X$ is a {\it holomorphic} line bundle over an $n$-dimensional compact complex manifold $X$, then its curvature form $\frac{i}{2\pi}\Theta_h(L)$ with respect to any $C^\infty$ Hermitian metric $h$ on $L$ is of type $(1,\,1)$. Hence, if $X$ is {\it projective}, we have: $$L \hspace{2ex} \mbox{is {\bf projectively divisorially nef}} \hspace{2ex} \iff c_1(L)^{n-1}.\{[D]\}_{BC}:=\int\limits_D\bigg(\frac{i}{2\pi}\Theta_h(L)\bigg)^{n-1}\geq 0$$ for all effective divisors $D\geq 0$ on $X$.

\vspace{2ex}

A well-known result (see e.g. [Ha70, $\S.6$, p. 34-36]) tells us that projectively divisorially nef holomorphic line bundles are, indeed, generalisations of nef such bundles.

\begin{The}\label{The:implication_nef_div-nef} Let $L\longrightarrow X$ be a holomorphic line bundle over a {\bf projective} manifold $X$. The following implication holds: \begin{eqnarray*}L \hspace{2ex} \mbox{is nef} \hspace{2ex} \implies \hspace{2ex} L \hspace{2ex} \mbox{is projectively divisorially nef}.\end{eqnarray*}

\end{The}

\noindent {\it Proof.} This follows at once from Kleiman's Theorem 6.1. in [Ha70, p. 34-36] which states that $L$ being nef is equivalent to $L^p.Y:=\int_Yc_1(L)^p\geq 0$ for every $p$-dimensional subvariety $Y\subset X$, for all $1\leq p\leq\mbox{dim}_\C X$. \hfill $\qed$

\section*{Acknowledgements} This work is part of the first-named author's PhD thesis under the supervision of the second-named author. The former wishes to express his gratitude to the latter for his constant guidance while this work was carried out, as well as to his Tunisian supervisor, Fathi Haggui, for constant support. Both authors are grateful to Luis Ugarte for helpful discussions about the Lie group $SL(2,\,\C)$ and related issues, as they are to Jean-Pierre Demailly for comments on an early version of the manuscript.

\vspace{3ex}

\noindent {\bf References.} \\

\noindent [Bro78]\, R. Brody --- {\it Compact Manifolds and Hyperbolicity} --- Trans. Amer. Math. Soc. {\bf 235} (1978), 213-219.  

\vspace{1ex}

\noindent [CY18]\, B.-L. Chen, X. Yang --- {\it Compact K\"ahler Manifolds Homotopic to Negatively Curved Riemannian Manifolds} --- Math. Ann. {\bf 370} (2018), DOI 10.1007/s00208-017-1521-7



\vspace{1ex}

\noindent [Dem92]\, J.-P. Demailly --- {\it Singular Hermitian Metrics on Positive Line Bundles} --- in Hulek K., Peternell T., Schneider M., Schreyer FO. (eds) ``Complex Algebraic Varieties''. Lecture Notes in Mathematics, vol {\bf 1507}, Springer, Berlin, Heidelberg.



\vspace{1ex}

\noindent [Dem00]\, J.-P. Demailly --- {\it Multiplier Ideal Sheaves and Analytic Methods in Algebraic Geometry} --- in ``School on Vanishing Theorems and Effective Results in Algebraic Geometry''; Trieste (Italy); 25 Apr -- 12 May 2000; p. 1-148.

\vspace{1ex}

\noindent [DP20]\, S. Dinew, D. Popovici --- {\it A Generalised Volume Invariant for Aeppli Cohomology Classes of Hermitian-Symplectic Metrics} --- Adv. Math. (2021), https://doi.org/10.1016/j.aim.2021.108056

\vspace{1ex}

\noindent [dTh10] H. De Th\'elin --- {\it Ahlfors' Currents in Higher Dimension} --- Ann. Fac. Sci. Toulouse {\bf XIX}, no. 1 (2010), p. 121-133. 

\vspace{1ex}

\noindent [FLY12]\, J. Fu, J. Li, S.-T. Yau -- {\it Balanced Metrics on Non-K\"ahler Calabi-Yau Threefolds} --- J. Differential Geom. 90 (2012), p. 81–129.

\vspace{1ex}

\noindent [Fri89]\, R. Friedman --{\it On Threefolds with Trivial Canonical Bundle} --- in Complex Geometry and Lie Theory (Sundance, UT, 1989), Proc. Sympos. Pure Math., vol. 53, Amer. Math. Soc., Providence, RI, 1991, p. 103–134.



\vspace{1ex}

\noindent [Gau77a]\, P. Gauduchon --- {\it Le th\'eor\`eme de l'excentricit\'e nulle} --- C. R. Acad. Sci. Paris, S\'er. A, {\bf 285} (1977), 387-390.

\vspace{1ex}

\noindent [Gau77b]\, P. Gauduchon --- {\it Fibr\'es hermitiens \`a endomorphisme de Ricci non n\'egatif} --- Bull. Soc. Math. France {\bf 105} (1977) 113-140.

\vspace{1ex}

\noindent [Gro91]\, M. Gromov --- {\it K\"ahler Hyperbolicity and $L^2$ Hodge Theory} --- J. Diff. Geom. {\bf 33} (1991), 263-292.

\vspace{1ex}

\noindent [Ha70]\, R. Hartshorne --- {\it Ample Subvarieties of Algebraic Varieties} --- Lecture Notes in Mathematics, no. {\bf 156}, Springer-Verlag, Berlin, 1970.

\vspace{1ex}

\noindent [Kob70]\, S. Kobayashi --- {\it Hyperbolic Manifolds and Holomorphic Mappings} --- Marcel Dekker, New York (1970). 

\vspace{1ex}

\noindent [Lam99]\, A. Lamari --- {\it Courants k\"ahl\'eriens et surfaces compactes} --- Ann. Inst. Fourier {\bf 49}, no. 1 (1999), 263-285.

\vspace{1ex}

\noindent [Lan87]\, S. Lang --- {\it Introduction to Complex Hyperbolic Spaces} --- Springer-Verlag (1987).

\vspace{1ex}

\noindent [LT93]\, P. Lu, G. Tian -- {\it The Complex Structures on Connected Sums of $S^3\times S^3$} --- in Manifolds and Geometry (Pisa, 1993), Sympos. Math., XXXVI, Cambridge Univ. Press, Cambridge, 1996, p. 284–293.

\vspace{1ex}

\noindent [Mic83]\, M. L. Michelsohn --- {\it On the Existence of Special Metrics in Complex Geometry} --- Acta Math. {\bf 143} (1983) 261-295.

\vspace{1ex}

\noindent [Nak75]\, I. Nakamura --- {\it Complex parallelisable manifolds and their small deformations} --- J. Differ. Geom. {\bf 10}, (1975), 85-112.

\vspace{1ex}

\noindent [Pop13]\, D. Popovici --- {\it Deformation Limits of Projective Manifolds: Hodge Numbers and Strongly Gauduchon Metrics} --- Invent. Math. {\bf 194} (2013), 515-534.

\vspace{1ex}

\noindent [Pop15a]\, D. Popovici --- {\it Aeppli Cohomology Classes Associated with Gauduchon Metrics on Compact Complex Manifolds} --- Bull. Soc. Math. France {\bf 143}, no. 4 (2015), p. 763-800.

\vspace{1ex}

\noindent [Pop15b] \, D. Popovici --- {\it Sufficient Bigness Criterion for Differences of Two Nef Classes} --- Math. Ann. {\bf 364} (2016), 649-655.

\vspace{1ex}

\noindent [PU18]\, D. Popovici, L. Ugarte --- {\it Compact Complex Manifolds with Small Gauduchon Cone} --- Proc. LMS {\bf 116}, no. 5 (2018) doi:10.1112/plms.12110.

\vspace{1ex}

\noindent [Sch07]\, M. Schweitzer --- {\it Autour de la cohomologie de Bott-Chern} --- arXiv e-print math.AG/0709.3528v1.

\vspace{1ex}

\noindent [ST10]\, J. Streets, G. Tian --- {\it A Parabolic Flow of Pluriclosed Metrics} --- Int. Math. Res. Notices, {\bf 16} (2010), 3101-3133.

\vspace{1ex}

\noindent [Voi02]\, C. Voisin --- {\it Hodge Theory and Complex Algebraic Geometry. I.} --- Cambridge Studies in Advanced Mathematics, 76, Cambridge University Press, Cambridge, 2002.

\vspace{1ex}

\noindent [Wan54]\, H.-C. Wang --- {\it Complex Parallisable Manifolds} --- 
Proc. Amer. Math. Soc. {\bf 5} (1954), 771--776.

\vspace{1ex}

\noindent [Yac98]\, A. Yachou --- {\it Sur les vari\'et\'es semi-k\"ahl\'eriennes} --- PhD Thesis, University of Lille.

\vspace{1ex}

\noindent [YZZ19]\, S.-T. Yau, Q. Zhao, F. Zheng --- {\it On Strominger K\"ahler-like Manifolds with Degenerate Torsion} --- arXiv e-print DG 1908.05322v2

\vspace{3ex}

\noindent Institut de Math\'ematiques de Toulouse, Universit\'e Paul Sabatier,

\noindent 118 route de Narbonne, 31062 Toulouse, France

\noindent Email: almarouanisamir@gmail.com AND popovici@math.univ-toulouse.fr

\vspace{2ex}

\noindent {\it For the first-named author only, also}:

\vspace{1ex}

\noindent Universit\'e de Monastir, Facult\'e des Sciences de Monastir

\noindent Laboratoire de recherche Analyse, G\'eom\'etrie et Applications LR/18/ES/16

\noindent Avenue de l'environnement 5019, Monastir, Tunisie

\end{document}